\documentclass{elsart}
\usepackage[dvips]{epsfig}
\newcommand{\R}{I \! \! R}
\newcommand{\N}{I \! \! {N}}
\newcommand{\C}{I \! \! \! \! {C}}

\renewcommand{\o}{\omega}
\newcommand{\uu}{{\underline u}}
\newcommand{\uy}{{\underline y}}
\newcommand{\uv}{{\underline v}}

\newcommand{\us}{{\underline s}}
\newcommand{\uc}{{\underline c}}

\newcommand{\unu}{{\underline \nu}}
\newcommand{\ux}{{\underline x}}

\newcommand{\uf}{{\underline f}}

\newcommand{\uaa}{{\underline a}}

\newcommand{\ue}{{\underline e}}

\newcommand{\uxi}{{\underline{\xi}}}
\newcommand{\uzet}{{\underline{\zeta}}}
\newcommand{\ugam}{{\underline{\gamma}}}
\newcommand{\ualf}{{\underline{\alpha}}}
\newcommand{\bb}{\begin{eqnarray}}
\newcommand{\be}{\end{eqnarray}}
\newcommand{\bbs}{\begin{eqnarray*}}
\newcommand{\bes}{\end{eqnarray*}}

\newtheorem{teo}{Theorem}

\newtheorem{lemma}{Lemma}

\newtheorem{defi}{Definition}

\begin{document}

\title{A new transform for solving the noisy complex
exponentials approximation problem}
\author{P. Barone }
\address{ Istituto per le Applicazioni del Calcolo ''M. Picone'',
C.N.R.,\\
Viale del Policlinico 137, 00161 Rome, Italy \\
e-mail: barone@iac.rm.cnr.it \\
fax: 39-6-4404306}

\maketitle

 \newpage

\section*{Abstract}

The problem of estimating a complex measure made up by a linear
combination of Dirac distributions centered on points of the complex
plane from a finite number of its complex moments affected by
additive i.i.d. Gaussian noise is considered. A random measure is
defined whose expectation approximates the unknown measure under
suitable conditions. An estimator of the approximating measure is
then proposed as well as a new discrete transform of the noisy
moments that allows to compute an estimate of the unknown measure. A
small simulation study is also performed to experimentally check the
goodness of the approximations.

{\it Key words and phrases:}
 Complex moments; Pade'  approximants; logarithmic potentials;
 random  determinants; random polynomials; pencils of matrices

\newpage

\section*{Introduction}

Let us consider the complex measure defined on a compact set
$D\subset\C$ by
$$S(z)=\sum_{j=1}^p c_j\delta(z-\xi_j),\;\;\xi_j\in \mbox{int}(
D), \;\;c_j\in\C$$  and let be
$$s_k=\int_Dz^kS(z)dz=\int\!\!\int_{\!\!\!\!\!\!D}(x+iy)^k S(x+iy)dx dy,
\;\;k=0,1,2,\dots$$ the complex moments. It turns out that
\begin{eqnarray}s_k=\sum_{j=1}^p c_j\xi_j^k.\label{modale}\end{eqnarray}
Let us assume to know an even number $n\ge 2p$ of noisy complex
moments
$${\bf a}_k=s_k+\mbox{\boldmath $\nu$}_k,\quad k=0,1,2,\dots,n-1
$$ where $\mbox{\boldmath $\nu$}_k$ is a complex Gaussian, zero mean, white noise,
 with finite known
variance $\sigma^2$. In the following all random quantities are
denoted by bold characters. We want to estimate $S(z)$ from $\{{\bf
a}_k\}_{k=0,\dots,n-1}$. From equation (\ref{modale}) this is
equivalent to estimate $p,c_j,\xi_j,\;j=1,\dots,p$, which is the
well known difficult problem of complex exponentials approximation.

The problem is central in many disciplines and appears in the
literature in different forms and contexts (see e.g.
\cite{donoho,gmv,osb,scharf,vpb}). The assumptions about the noise
variance (constant and known) are made here to simplify the
analysis. However in many applications the noise is an instrumental
one which is well represented by a white noise, zero mean, Gaussian
process whose variance is known or easy to estimate. A typical
example is provided by NMR spectroscopy (see e.g. \cite{farr}).

In the noiseless case the problem becomes the complex exponential
interpolation problem \cite{hen2}. Conditions for existence and
unicity of the solution are (\cite[Th.7.2c]{hen2}): $$det
U_0(\us)\ne 0,\;\;\; det U_1(\us) \ne 0$$  where
$$U(s_0,\dots,s_{2p-2})=\left[\begin{array}{llll}
s_0 & s_{1} &\dots &s_{p-1} \\
s_{1} & s_{2} &\dots &s_{p} \\
. & . &\dots &. \\
s_{p-1} & s_{p} &\dots &s_{2p-2}
  \end{array}\right]$$ and $$U_0(\us)=U(s_0,\dots,s_{2p-2}),\;\;\;\;U_1(\us)=U(s_1,\dots,s_{2p-1}).$$
In fact exactly $n=2p$ noiseless moments are sufficient to fully
retrieve $S(z)$, where $$p=\max_{n\in\N}\{n\;
|\;det(U(s_0,\dots,s_{n-2}))\neq 0\}.$$ Moreover
$(\xi_j,\;j=1,\dots,p)$ are the generalized eigenvalues of the
pencil $P=[U_1(\us),U_0(\us)]$ i.e. they are the roots of the
polynomial in the variable $z$
$$det[U_1(\us)-zU_0(\us)]$$
 and $c_j$ are related to the generalized
eigenvector $\uu_j$ of $P$ by $c_j=\uu_j^T[s_0,\dots,s_{p-1}]^T$. In
fact from equation (\ref{modale}) we have
$\uc=V^{-1}[s_0,\dots,s_{p-1}]^T$ where
$$V=Vander(\xi_1,\dots,\xi_p)$$ is the square
Vandermonde matrix based on $(\xi_1,\dots,\xi_p)$. But it easy to
show (see e.g. \cite{bama2}) that
$$U_0(\us)=VCV^T,\;\;
U_1(\us)=VCZV^T$$ where
$$C=diag\{c_1,\dots,c_p\} \mbox{ and }
Z=diag\{\xi_1,\dots,\xi_p\}.$$ Therefore $\uu_k=V^{-T}\ue_k$ is the
right generalized eigenvector of $P$ corresponding to $\xi_k$, where
$\ue_k$ is the $k-$th column of the identity matrix $I_{p}$ of order
$p$.

 Viceversa when $s_k=0,\;\forall
k$ it was proved in \cite{cinese} that
$$det[U({\bf a}_0,\dots,{\bf a}_{n-2})]=det[U_0({\bf \uaa})]\neq 0\;\;\forall n \mbox{
a.s.}$$ and
$$det[U({\bf a}_1,\dots,{\bf a}_{n-1})]=det[U_1({\bf \uaa})]\neq 0\;\;\forall n \mbox{
a.s.}.$$ Moreover associated to the random polynomial \bb
det[U_1({\bf \uaa})-zU_0({\bf \uaa})]\label{rapol}\be a condensed
density $h_n(z)$ can be considered which is the expected value of
the (random) normalized counting measure on the zeros of this
polynomial i.e.
$$h_n(z)= \frac{2}{n}E\left[\sum_{j=1}^{n/2}\delta(z-\mbox{\boldmath $\xi$}_j)\right].
$$ It was proved in \cite{barja} that if $z=r e^{i\theta}$,  the
marginal condensed density $h_n^{(r)}(r)$ w.r. to $r$  of the
generalized eigenvalues is asymptotically in $n$ a Dirac $\delta$
supported on the unit circle $\forall \sigma^2$. Moreover for finite
$n$ the the marginal condensed density w.r. to $\theta$  is
uniformly distributed on $[-\pi,\pi]$. Starting from the generalized
eigenvalues $\mbox{\boldmath $\xi$}_j$ and generalized eigenvectors
${\bf \uu}_j$ of the pencil
$${\bf
P}=[U({\bf a}_1,\dots,{\bf a}_{n-1}),U({\bf a}_0,\dots,{\bf
a}_{n-2})]$$ we then define a family of random measures
$${\bf S}_n(z)=\sum_{j=1}^{n/2}{\bf c}_j\delta(z-\mbox{\boldmath $\xi$}_j)$$
where ${\bf c}_j={\bf \uu}_j^T[{\bf a}_0,\dots,{\bf a}_{n/2-1}]^T$
and we give conditions under which $E[{\bf S}_n(z)]$ approximates
$S(z)$. Moreover we define a
 discrete transform ($\mathsf{P}$-Transform)  on a lattice of points on $D$,
  which is an unbiased and consistent estimator of
   $E[{\bf S}_n(z)]$ on the lattice thus providing a
computational device to solve the original problem.

In \cite{baram} the same problem was afforded. The joint
distribution of the coefficients of the random polynomial
(\ref{rapol}) (when $s_k\ne 0,\;\forall k$) was approximated by  a
multivariate Gaussian distribution and a theorem by Hammersley
\cite{ham} was used to compute the associated condensed density of
its roots. An heuristic algorithm was then used to identify the main
peaks of the condensed density and to get estimates of $p,\;\xi_j
\mbox{ and } c_j,\;j=1,\dots,p$ based on them. In the present work
the ideas presented in \cite{baram} are put on a more rigorous
mathematical framework. A different approximation of the condensed
density is considered and an automatic estimation procedure is
proposed.

The paper is organized as follows. In the first section we study the
distribution of the generalized eigenvalues of the random pencil
${\bf P}$ and we give an easily computable approximate expression of
the associated condensed density. In section 2 we consider the
identifiability problem for $S(z)$ given the data ${\bf \uaa}$.
Conditions for identifiability are given and the approximation
properties of $E[{\bf S}_n(z)]$ are proved. In section 3 the
$\mathsf{P}$-transform is defined and its statistical properties are
studied. In section 4 the  procedure for estimating the parameters
$p$, $\{\xi_j, c_j,j=1,\dots,p\}$ of the unknown measure from the
$\mathsf{P}$-transform is described. Finally in section 5 some
experimental results on synthetic data are reported.

\section{Distribution of the generalized eigenvalues of the pencil ${\bf P}$}

We start by   making some technical assumptions on the noise model.
When $s_k=0\;\forall k$, we noticed in the introduction that ${\bf
\xi}_j$ are, asymptotically on $n$, uniformly distributed on the
unit circle. Therefore, when $s_k\ne 0$ is given by (\ref{modale}),
we can assume that  $n_p=n/2-p$ among the ${\bf
\xi}_j,j=1,\dots,n/2$ are related to noise  and then they can be
modeled for large $n$ by $\tilde{\xi}_j=e^{\frac{2\pi i j}{n_p}}$
i.e. by  uniformly spaced deterministic generalized eigenvalues.
Therefore the Vandermonde matrix based on
$\tilde{\xi}_j,j=1,\dots,n_p$ is simply given by $V=\sqrt{n_p}\cdot
F\in\C^{n_p \times n_p}$ where
$F_{hk}=\frac{1}{\sqrt{n_p}}e^{\frac{2\pi i hk}{n_p}}$ is the
discrete Fourier transform matrix. Hence
$${\bf \tilde{\uc}}= V^{-1}[\mbox{\boldmath $\nu$}_0,\dots,\mbox{\boldmath $\nu$}_{n_p-1}]^T=
\frac{1}{\sqrt{n_p}}F^H[\mbox{\boldmath
$\nu$}_0,\dots,\mbox{\boldmath $\nu$}_{n_p-1}]^T$$ and ${\bf
\tilde{\uc}}$ has a complex multivariate Gaussian distribution with
$$E[\tilde{{\bf c}}_j]=0 \mbox{ and }
E[\tilde{{\bf c}}_j\overline{\tilde{{\bf c}}_h}]=
\frac{\sigma^2}{n_p}\delta_{jh}.$$ Based on these observations we
define a new noise process as
$$\mbox{\boldmath $\tilde{\nu}$}_k=\left\{\begin{array}{ll}
\sum_{j=1}^{n_p} \tilde{{\bf c}}_j\tilde{\xi}_j^k, & \;k<n_p \\
\mbox{\boldmath $\nu$}_k, & \;k\ge n_p \end{array} \right. $$
 and we assume that ${\bf \tilde{\uc}}$ is independent of $\mbox{\boldmath $\nu$}_k,k\ge n_p.$
 But then
$E[\mbox{\boldmath $\tilde{\nu}$}_k]=0$ and
$$E[\mbox{\boldmath $\tilde{\nu}$}_k\overline{\mbox{\boldmath $\tilde{\nu}$}}_h]=\left\{\begin{array}{ll}
\sum_{i,j}^{1,n_p}\tilde{\xi}_i^k\overline{\tilde{\xi}}_j^hE[\tilde{{\bf
c}}_i\overline{\tilde{{\bf c}}}_j]=
\frac{\sigma^2}{n_p}\sum_{r=1}^{n_p}e^{\frac{2\pi i r(k-h)}{n_p}}=\sigma^2\delta_{hk}, & \;k,h<n_p \\
\sum_{j=1}^{n_p} E[\tilde{{\bf c}}_j\overline{\mbox{\boldmath $\nu$}}_h]\tilde{\xi}_j^k=0, & \;h\ge n_p,k< n_p \\
E[\mbox{\boldmath $\nu$}_k\overline{\mbox{\boldmath
$\nu$}}_h]=\sigma^2\delta_{hk}, & h,k \ge n_p
\end{array} \right. $$
We have then proved the following
\begin{lemma}
The random vectors $\mbox{\boldmath $\nu$}_k$ and $\mbox{\boldmath
$\tilde{\nu}$}_k,\;\;k=0,\dots,n-1$ are equal in
distribution.\label{lem1}
\end{lemma}
As a consequence  in the following we will use $\mbox{\boldmath
$\tilde{\nu}$}_k$ without loss of generality.

\noindent \begin{rem} We notice that when $s_k\ne 0$, if the
signal-to-noise ratio is defined as $SNR=
\frac{1}{\sigma}\min_{h=1,p}|c_h|$ we have
$$E[|\tilde{{\bf c}}_j|^2]=\frac{\sigma^2}{n_p}=\frac{\min_{h=1,p}|c_h|^2}
{n_p SNR^2}.$$  If $SNR\gg\sqrt{\frac{1}{n_p}}$ then $E[|\tilde{{\bf
c}}_j|^2]\ll |c_k|^2,\;\;\forall j,k. $
\end{rem}

A basic result which will be used extensively in the following is
given by
\begin{lemma}
Let $T=(T^{(1)},T^{(2)})$ be the transformation that maps every
realization $\uaa(\o)$ of ${\bf \uaa}$ to $(\uxi(\o), \uc(\o))$
given by $
a_k(\o)=\sum_{j=1}^{n/2}c_j(\o)\xi_j(\o)^{k},\;\;k=0,\dots,n-1,$
where $\o\in\Omega$ and $\Omega$  is the space of events. Then ${\bf
T}$ is a.s. one-to-one. Moreover,  for $\sigma\rightarrow 0$ and for
$j=1,\dots,n/2$
$$E[\mbox{\boldmath $\xi$}_j]=\left\{\begin{array}{ll}\xi_j+o(\sigma) & j=1,\dots,p \\
\tilde{\xi}_{j-p}+o(\sigma),\;\; & j=p+1,\dots,n/2\end{array}\right
. $$
$$E[{\bf c}_j]=\left\{\begin{array}{ll}c_j +o(\sigma),& j=1,\dots,p \\
o(\sigma),\;\; & j=p+1,\dots,n/2\end{array}\right . $$ \label{lem11}
\end{lemma}

\noindent \underline{proof}

\noindent From \cite{cinese} we know that a.s.
$det[U_h(\mbox{\boldmath $\unu$})]\neq 0,\;h=0,1$. Moreover, with
probability $1$, there is no functional dependence between
$\mbox{\boldmath $\unu$}$ and $\us$. Therefore a.s. $det[U_h({\bf
\uaa})]\neq 0,\;h=0,1$. But then a.s. the complex exponential
interpolation problem for ${\bf \uaa}$ has an unique solution
$\forall \o$ hence ${\bf T}$ is a.s. one-to-one. The second part of
the thesis is based on a Taylor expansion of ${\bf T}$ around a
suitable point $\ux_0$. A natural candidate for $\ux_0$ would be
$\us$. However we notice that $T^{(1)}(\us)$ is not defined if
$n>2p$, and, as a consequence, also $T^{(2)}(\us)$ is not defined in
this case. Therefore, by using Lemma \ref{lem1}, without loss of
generality, we assume that the noise is represented by
$\mbox{\boldmath $\tilde{\nu}$}_k$ i.e.
$${\bf a}_k=\left\{\begin{array}{ll}\sum_{j=1}^p
c_j\xi_j^k+\sum_{j=p+1}^{n/2}
\tilde{{\bf c}}_{j-p}\tilde{\xi}_{j-p}^k ,\;\; & k=0,\dots,n_p-1 \\
\sum_{j=1}^p c_j\xi_j^k+\mbox{\boldmath $\nu$}_k,\;\; &
k=n_p,\dots,n-1\end{array}\right.$$ where  $n_p=n/2-p.$ We  then
define a new sequence $\tilde{s}_k$ by
$$\tilde{s}_k=\sum_{j=1}^p c_j\xi_j^k+\sigma^\alpha\sum_{j=p+1}^{n/2} \tilde{\xi}_{j-p}^k ,\;\;
\alpha\ge 2,\;\;  k=0,\dots,n-1$$
 and we
consider the process ${\bf a}_k$  as a perturbation of
 $\tilde{s}_k$. Therefore we choose $\ux_0=\tilde{\us}$
and notice that
$$T^{(1)}(\tilde{\us})_j=\left\{\begin{array}{ll}\xi_j & j=1,\dots,p \\
\tilde{\xi}_{j-p},\;\; & j=p+1,\dots,n/2\end{array}\right . $$
$$T^{(2)}(\tilde{\us})_j=\left\{\begin{array}{ll}c_j & j=1,\dots,p \\
\sigma^\alpha,\;\; & j=p+1,\dots,n/2\end{array}\right . $$ \noindent
We now prove that each component of $T^{(1)}(\uaa)$ is an analytic
function of $\uaa$ when  $\uaa$ belong to small neighbor of
$\tilde{\us}$. The proof follows closely \cite{sb}[Th.6.9.8]. For
each fixed $\o$, the polynomial
$$\phi(z,\uaa)=det[U_1(\uaa)-zU_0(\uaa)]$$
is an analytic function of $z$ and $\uaa$. Let $\hat{\xi}$ be a zero
of $\phi(z,\tilde{\us})$ and
$$K=\{\zeta||\zeta-\hat{\xi}|=r\},\;\;r>0$$ be a circle around $\hat{\xi}$ not
containing any other generalized eigenvalue of the pencil
$$\tilde{
P}=[U(\tilde{s}_1,\dots,\tilde{s}_{n-1}),U(\tilde{s}_0,
\dots,\tilde{s}_{n-2})].$$ We want to show that $K$ does not pass
through any zero of $\phi(z,\uaa)$. In fact by the definition of $K$
it follows that
$$\inf_{\zeta\in K}|\phi(\zeta,\tilde{\us})|>0.$$ But $\phi(z,\uaa)$ depends
continuously on $\uaa$, hence there exists
$B=\{\ux\in\C^n||\ux-\tilde{\us}|<\rho\},\;\;\rho>0$ such that
$$\inf_{\zeta\in K}|\phi(\zeta,\uaa)|>0,\;\;\forall \uaa\in B.$$ By the
principle of argument, the number of zeros of $\phi(z,\uaa)$ within
$K$ is given by $$N(\uaa)=\frac{1}{2\pi
i}\oint_K\frac{\phi'(z,\uaa)}{\phi(z,\uaa)}dz,\;\;\phi'=\frac{\partial\phi}{\partial
z}$$ which is continuous in $B$; hence $$1=N(\tilde{\us})=N(\uaa),
\;\;\uaa\in B.$$ Moreover the simple zero $\xi(\o)$ of
$\phi(z,\uaa)$ inside $K$ admits the representation (see e.g.
\cite{mark})
$$\xi(\o)=\frac{1}{2\pi
i}\oint_K\frac{z\phi'(z,\uaa)}{\phi(z,\uaa)}dz .$$ For $\uaa\in B$
the integrand is an analytic function of $\uaa$ and therefore also
$\xi(\o)$ is an analytic function of $\uaa$ when $\uaa\in B$.

\noindent We now consider $T^{(2)}(\uaa)$. We notice that each
component  can be obtained as a rational function of the components
of $T^{(1)}(\uaa)$ by the formula $ c_j=\ue_j^TV^{-H}
\uaa,\;j=1,\dots,n/2$ where $V$ is the Vandermonde matrix based on
$T^{(1)}(\uaa)$. Therefore also $ c_j$ is an analytic function of
$\uaa$ when $\uaa\in B$.

\noindent As $T^{(h)}=T_R^{(h)}+iT_I^{(h)}$ is analytic for $\uaa\in
B$, $T_R^{(h)}$ and $T_I^{(h)}$ are real analytic functions of
$\uaa_R,\uaa_I$ where $\uaa=\uaa_R+i\uaa_I$, (e.g.
\cite{gk}[pg.99]). Therefore they admit a Taylor series expansion
 around $\tilde{\us}$ when $\uaa\in B$: \bb
T_{Rk}^{(h)}(\uaa)&=&T_{Rk}^{(h)}(\tilde{\us}) + \sum_{i=0}^{n-1}
\frac{\partial T_{Rk}^{(h)}(\uaa) }{\partial
a_{Ri}}_{|\uaa=\tilde{\us}}[a_{Ri}-\tilde{s}_{Ri}]+ \nonumber
\\& & \sum_{i=0}^{n-1}
\frac{\partial T_{Rk}^{(h)}(\uaa) }{\partial
a_{Ii}}_{|\uaa=\tilde{\us}}[a_{Ii}-\tilde{s}_{Ii}] + \nonumber
\\& &
\frac{1}{2}\sum_{i=0}^{n-1}\sum_{j=0}^{n-1} \frac{\partial^2
T_{Rk}^{(h)}(\uaa) }{\partial a_{Ri}\partial
a_{Rj}}_{|\uaa=\tilde{\us}}[a_{Ri}-\tilde{s}_{Ri}][a_{Rj}-\tilde{s}_{Rj}]+
\nonumber
\\& &
\frac{1}{2}\sum_{i=0}^{n-1}\sum_{j=0}^{n-1} \frac{\partial^2
T_{Rk}^{(h)}(\uaa) }{\partial a_{Ii}\partial
a_{Ij}}_{|\uaa=\tilde{\us}}[a_{Ii}-\tilde{s}_{Ii}][a_{Ij}-\tilde{s}_{Ij}]+
\nonumber
\\& &
\sum_{i=0}^{n-1}\sum_{j=0}^{n-1} \frac{\partial^2 T_{Rk}^{(h)}(
\uaa) }{\partial a_{Ri}\partial
a_{Ij}}_{|\uaa=\tilde{\us}}[a_{Ri}-\tilde{s}_{Ri}][a_{Ij}-\tilde{s}_{Ij}]+...\nonumber
\be and analogously for $T_{Ik}^{(h)}(\uaa)$. Taking expectations we
get  \bb & &E\left[({\bf a}_{Ri}-\tilde{s}_{Ri})\right]
 = \left[s_{Ri}-\tilde{s}_{Ri}\right]=\sigma^\alpha  \cdot
C_i,\;\;C_i=\sum_{j=p+1}^{n/2}\tilde{\xi}^i_{j-p}\nonumber\be \bb &
&E\left[({\bf a}_{Ri}-\tilde{s}_{Ri})({\bf
a}_{Rj}-\tilde{s}_{Rj})\right]= E\left[({\bf
a}_{Ri}-s_{Ri}+\sigma^\alpha C_i)({\bf a}_{Rj}-s_{Rj}+\sigma^\alpha
C_j)\right] \nonumber
\\&=&
\frac{\sigma^2}{2}\delta_{ij}+\sigma^{2\alpha}C_iC_j\nonumber\be and
analogously for the other terms. Remembering the independence of the
real and imaginary parts of ${\bf a}_k$, we finally get
$$E[T_k^{(h)}({\bf \uaa})]=T_k^{(h)}(\tilde{\us})+o(\sigma).\qed$$

We start now the study of the distribution in $\C$ of the
generalized eigenvalues of ${\bf P}$ by making some qualitative
statements already present in the literature.
 For each realization $\o$, let
$\{c_j(\o),\xi_j(\o)\},\;j=1,\dots,n/2$ be the solution of  the
complex exponential interpolation problem for the data
$a_k(\o),\;k=0,\dots,n-1$. It is well known that we can then define
the  Pade' approximant $$[n/2-1,n/2](z,\o)=z\sum_{j=1}^{n/2}
\frac{c_j(\o)}{z-\xi_j(\o)}=Q_{n/2-1}(z^{-1})/P_{n/2}(z^{-1})$$ to
the $Z-$transform of $\{a_k(\o)\}$ given by
$$f(z,\o)=\sum_{k=0}^{\infty}a_k(\o)z^{-k}=f_s(z)+f_{\nu}(z,\o)$$ where
$$f_s(z)=\sum_{k=0}^{\infty}s_kz^{-k}= \sum_{j=1}^{p}
c_j\sum_{k=0}^{\infty}(\xi_j/z)^k= z\sum_{j=1}^{p}
\frac{c_j}{z-\xi_j},\;\;|z|>1$$ and, because of Lemma \ref{lem1},
$$f_{\nu}(z,\o)\approx z\sum_{j=1}^{n_p}
\frac{\tilde{c}_j(\o)}{z-\tilde{\xi}_j}$$ $f(z,\o)$ is then defined
outside the unit circle and can be extended to $D$ by analytic
continuation.  We get then
$$f(z,\o)\approx z\tilde{q}_{n/2-1}(z)/\tilde{p}_{n/2}(z)=
\frac{z\prod_{j=1}^{n/2-1}(z-\delta_j(\o))}{
\prod_{j=1}^{p}(z-\xi_j)\prod_{j=1}^{n_p}(z-\tilde{\xi}_j)}$$ and
$$g(z,\o)=\log(z^{-1}f(z,\o))=\sum_{j=1}^{n/2-1}\log(z-\delta_j(\o))-
\sum_{j=1}^{p}\log(z-\xi_j)-\sum_{j=1}^{n_p}
\log(z-\tilde{\xi}_j).$$

We want to study the location in $\C$ of $\xi_j(\o)$. To this aim,
following \cite{maba2}, we remember that $p_n(z)=z^nP_{n}(z^{-1})$
satisfy the following orthogonality relation
$$\int_\Gamma z^{-1}f(z,\o)p_n(z)z^kdz=0,\;\;k=0,\dots,n-1$$
where $\Gamma$ is a union of closed curves enclosing the poles of
$f(z,\o)$ i.e. the numbers $\xi_j,\;j=1\dots,p$ and
$\tilde{\xi}_j,\;\;j=1,\dots,n_p$. By using the Szego integral
representation of such polynomials and a saddle point argument, it
turns out that the Pade' poles $\xi_j(\o),\;j=1,\dots,n/2$ ,
asymptotically on $n$, satisfy the following system of algebraic
equations \begin{eqnarray*}2\sum_{j\ne
k}^{1,n/2}\frac{1}{(\xi_k(\o)-\xi_j(\o))}+g'(\xi_k(\o))=0
\;\;\;\;\;\;\;\;\;\;\;\;\;k=1,\dots,n/2\end{eqnarray*} or
\begin{eqnarray*}2\sum_{j\ne
k}^{1,n/2}\frac{1}{(\xi_k(\o)-\xi_j(\o))}
+\sum_{j=1}^{n/2-1}\frac{1}{(\xi_k(\o)-\delta_j(\o))}+\\-
\sum_{j=1}^{p}\frac{1}{(\xi_k(\o)-\xi_j)}-\sum_{j=1}^{n_p}
\frac{1}{(\xi_k(\o)-\tilde{\xi}_j)}=0,\;\;k=1,\dots,n/2
\end{eqnarray*}
These equations can be interpreted as conditions of electrostatic
equilibrium of a set of  charges in the presence of an electric
external field corresponding to $g'(z,\o)$. Therefore the Pade'
poles $\xi_k(\o)$ are attracted by $\xi_j,\;j=1,\dots,p \mbox{ and
}\tilde{\xi}_j,\;j=1,\dots n_p$  and they are repelled by each other
and by the  zeros $\delta_j(\o)$ of $\tilde{q}_{n/2-1}(z)$. However
\bb\tilde{q}_{n/2-1}(z)&=&\sum_{j=1}^p c_j\prod_{k\ne
j}^{1,p}(z-\xi_k)
\prod_{k=1}^{n_p}(z-\tilde{\xi}_k)\\&+&\sum_{j=1}^{n_p}
\tilde{c}_j(\o)\prod_{k=1}^{p}(z-\xi_k) \prod_{k\ne
j}^{1,n_p}(z-\tilde{\xi}_k).\be As $\forall
\o,\;\;|\tilde{c}_j(\o)|^2\ll \min_h |c_h|^2$ if the SNR is
sufficiently high (see Remark after Lemma \ref{lem1}), we can
approximate $\tilde{q}_{n/2-1}(z)$ by
$$\prod_{k=1}^{n_p}(z-\tilde{\xi}_k)\sum_{j=1}^p c_j\prod_{k\ne j}^{1,p}(z-\xi_k)
$$
 hence  $n_p$
zeros are close to $\tilde{\xi}_k$,  and the other $p-1$ are close
to the zeros of the polynomial
$$q_{p-1}(z)=\sum_{j=1}^pc_j\prod_{k\ne j}^{1,p}(z-\xi_k)$$
which is the numerator of $z^{-1}f_s(z)$. We notice that if
$|c_h|\ll |c_k|,\;\forall k\ne h$ then
$$q_{p-1}(z)\approx \sum_{j\ne h}^{1,p}c_j\prod_{k\ne j}^{1,p}(z-\xi_k)
=(z-\xi_h)\sum_{j=1}^pc_j\prod_{k\ne j,h}^{1,p}(z-\xi_k)$$ Hence,
because of the continuous dependence of the roots from the
coefficient of a polynomial, $q_{p-1}(z)$ has a zero as close to
$\xi_h$ as $|c_h|$ is small with respect to $|c_k|,\; k\ne h$.
 Therefore the Pade' poles $\xi_k(\o)$
\begin{itemize}
\item
are attracted by $\xi_j,\;j=1,\dots,p $
\item are attracted
by $\tilde{\xi}_j,\;j=1,\dots n_p$
\item are repelled from $\xi_j(\o),\;j\ne k$
\item are repelled from $\tilde{\xi}_j,\;j=1,\dots n_p$

\item are repelled from other $p-1$ points in the complex
plane which are as close to $\xi_j$ as $|c_j|$ is small with respect
to $|c_h|,\; h\ne j$.
\end{itemize}

 Summing up
a $\xi_k$ with a large $|c_k|$ will attract a Pade' pole without
being disturbed by the repulsion exerted by the zeros of
$\tilde{q}_{n/2-1}(z)$. Moreover close to such a point a gap of
Pade' poles can be expected because of the repulsion exerted by
Pade' poles to each other. A $\xi_k$ with a small $|c_k|$ will still
attract a Pade' pole but not so close because of the repulsion
exerted by a close zero.  The Pade' poles not related to the signal
are expected to be attracted by $\tilde{\xi}_k$ which at the same
time will repel them. Moreover they are repelled by $\xi_k$ hence
they are likely to be located in between $\tilde{\xi}_k$ and far
from $\xi_k$. A picture of this behavior is given in fig.1. We
notice that the qualitative results discussed above are consistent
with those obtained in \cite{bama1} under a more stringent
hypothesis about the noise.

We now wish to define a mathematical tool to quantify these
qualitative statements. To this aim we remember that
$\mbox{\boldmath $\xi$}_k,k=1,\dots,n/2$ are the generalized
eigenvalues of the pencil ${\bf P}$ and therefore they satisfy the
equation
$${\bf P}_{n/2}(z^{-1})=det[U_1({\bf \uaa})-zU_0({\bf \uaa})]=0.$$
Then a condensed density $h_n(z)$ can  be considered which is the
expected value of the (random) normalized counting measure on the
zeros of this polynomial i.e.
$$h_n(z)= \frac{2}{n}E\left[\sum_{j=1}^{n/2}\delta(z-\mbox{\boldmath $\xi$}_j)\right].
$$
The following theorem holds whose proof is the same of that of
Theorem 1 in \cite{barja}:
\begin{teo}
The condensed density of the zeros of the random polynomial ${\bf
Q}(z)={\bf P}_{n/2}(z^{-1})$ is given by
\begin{eqnarray}h_n(z)=\frac{1}{4\pi}\Delta  u_n(z)
\label{eqf}\end{eqnarray}  where $\Delta$ denotes the Laplacian
operator  with respect to $x,y$ if $z=x+iy$ and
\begin{eqnarray}u_n(z)=\frac{2}{n}E\left\{\log(|{\bf Q}(z)|^2)\right\}\label{un}\end{eqnarray}
\end{teo}
The condensed density provides the required quantitative information
about the distribution of the Pade' poles in the complex plane. If
the SNR is sufficiently high, after the qualitative statements made
above about the location of the Pade' poles,  a peak of $h_n(z)$ can
be expected in a neighborhood of each of the complex exponentials
$\xi_k,k=1,\dots,p$ and the volume under the peak gives the
probability of finding a Pade' pole in that neighborhood. This is
confirmed by the following
\begin{teo}
If $\sigma>0$, the condensed density $h_n(z, \sigma)$ is a
continuous function of $z$  given by \bb
h_n(z,\sigma)=\frac{2}{n(\pi\sigma^2)^{n}}\sum_{j=1}^{n/2}
\int_{\C^{n/2-1}}\int_{\C^{n/2}}
J_C^*(\uzet^*_j,z,\ugam)e^{-\frac{1}{\sigma^2}
\sum_{k=0}^{n-1}|\sum_{h\ne j}^{1,n/2} \gamma_h\zeta_h^{k}+\gamma_j
z^{k}-s_k|^2} d\uzet^*_jd\ugam \label{integr} \be where
 $\uzet^*_j=\{\zeta_h,h\ne j\}$ and $$J_C^*(\uzet^*_j,z,\ugam)=
\left\{\begin{array}{cc} \gamma & \mbox{ if } n=2\\
(-1)^{n/2}\prod_{ j=1}^{1,n/2}\gamma_j\prod_{r<h,r\ne
j}(\zeta_r-\zeta_h)^4\prod_{r\ne j}(\zeta_r-z)^4 & \mbox{ if } n\ge
4\end{array}\right.$$  Moreover
 $h_{n}(z, \sigma)$ converges weakly to the positive measure
$\frac{2}{n}\sum_{j=1}^p\delta(z-\xi_j)\mbox{ when
}\sigma\rightarrow 0$. \label{teo2}
\end{teo}
\noindent\underline{proof}

Let us consider the transformation
$T_n:\ualf\rightarrow(\uzet,\ugam)$ given by
$$\alpha_k=\sum_{j=1}^{n/2}\gamma_j\zeta_j^{k}$$
or
$$(T_n^{(1)}(\ualf))_j=\zeta_j,\;\;(T_n^{(2)}(\ualf))_j=\gamma_j.$$
In the following, to simplify notations, $(T_n^{(1)}(\ualf))_j$ will
be denoted by $\zeta_j(\ualf)$. We have \bb
h_n(z,\sigma)&=&\frac{2}{n}E\left[\sum_{j=1}^{n/2}\delta(z-\mbox{\boldmath
$\xi$}_j)\right]\\&=&
\frac{2}{n(\pi\sigma^2)^{n}}\sum_{j=1}^{n/2}\displaystyle
\int_{\C^n}\delta(z-\zeta_j(\ualf)
e^{-\frac{1}{\sigma^2}\sum_{k=0}^{n-1}|\alpha_k-s_k|^2}d\ualf;
\label{dc1}\be As the complex Jacobian of $T_n^{-1}$ is (see
\cite{flha,krat}) ($n$ was assumed even):
$$J_C(\uzet,\ugam)=\left\{\begin{array}{cc} \gamma & \mbox{ if } n=2\\
 (-1)^{n/2}\prod_{j=1}^{n/2}\gamma_j
\prod_{j<h}(\zeta_j-\zeta_h)^4 & \mbox{ if } n\ge
4\end{array}\right.,$$ by making a change of variables we have
\begin{eqnarray*}h_n(z,\sigma)&=&\frac{2}{n(\pi\sigma^2)^{n}}
\sum_{j=1}^{n/2}\int_{\C^{n/2}}\int_{\C^{n/2}}\delta(z-\zeta_j)
J_C(\uzet,\ugam)e^{-\frac{1}{\sigma^2}\sum_{k=0}^{n-1}|
\sum_{h=1}^{n/2}\gamma_h\zeta_h^{k}-s_k|^2}d\uzet
d\ugam\\
&=&\frac{2}{n(\pi\sigma^2)^{n}}\sum_{j=1}^{n/2}
\int_{\C^{n/2-1}}\int_{\C^{n/2}}
J_C^*(\uzet^*_j,z,\ugam)e^{-\frac{1}{\sigma^2}
\sum_{k=0}^{n-1}|\sum_{h\ne j}^{1,n/2}
\gamma_h\zeta_h^{k}+\gamma_j
z^{k}-s_k|^2} d\uzet^*_jd\ugam
\end{eqnarray*} where
 $\uzet^*_j=\{\zeta_h,h\ne j\}$ and $$J_C^*(\uzet^*_j,z,\ugam)=
\left\{\begin{array}{cc} \gamma & \mbox{ if } n=2\\
(-1)^{n/2}\prod_{j=1}^{n/2}\gamma_j\prod_{r<h,r\ne
j}(\zeta_r-\zeta_h)^4\prod_{r\ne j}(\zeta_r-z)^4 & \mbox{ if } n\ge
4\end{array}\right.$$ The integral above
 converges uniformly for $z\in D$, hence $h_n(z)$ is continuous in
 $D$.
We prove now  that $h_{2p}(z, \sigma)$ converges weakly to
$\frac{1}{p}\sum_{j=1}^p\delta(z-\xi_j)$ when $\sigma\rightarrow 0$.
Let $\Phi(z)\in C^\infty$ be a bounded test function supported on
$\C$. We have
\begin{eqnarray*}
& & \int_{\C} h_{2p}(z,\sigma)\Phi(z)dz \\&=&
\frac{1}{p(\pi\sigma^2)^{2p}}\sum_{j=1}^{p}\displaystyle \int_{\C}
\Phi(z)\left[\int_{\C^{2p}}\delta(z-\zeta_j(\ualf))
e^{-\frac{1}{\sigma^2}\sum_{k=0}^{2p-1}|\alpha_k-s_k|^2}d\ualf\right]
dz\\& =&
 \frac{1}{p(\pi\sigma^2)^{2p}}\sum_{j=1}^{p}\displaystyle
\int_{\C^{2p}}\Phi(\zeta_j(\ualf))
e^{-\frac{1}{\sigma^2}\sum_{k=0}^{2p-1}|\alpha_k-s_k|^2}d\ualf \\&=&
 \frac{1}{p}\sum_{j=1}^{p}\displaystyle
\int_{\C^{2p}}\Phi(\zeta_j(\uy\sigma+\us))
\frac{e^{-\sum_{k=0}^{2p-1}|\uy_k|^2}}{\pi^{2p}}d\uy.
\end{eqnarray*}
As  $\Phi(z)$ is continuous and bounded and $\zeta_j$ is analytic in
a neighbor of $\us$ by Lemma \ref{lem11}, by the dominated
convergence theorem we get
\begin{eqnarray*}
\lim_{\sigma\rightarrow 0}\int_\Omega h_{2p}(z,\sigma)\Phi(z)dz=
 \frac{1}{p}\sum_{j=1}^{p}\displaystyle
\int_{\C^{2p}}\lim_{\sigma\rightarrow 0}\Phi(\zeta_j(\uy\sigma+\us))
\frac{e^{-\sum_{k=0}^{2p-1}|\uy_k|^2}}{\pi^{2p}}d\uy =\\
\frac{1}{p}\sum_{j=1}^{p}\Phi(\zeta_j(\us))\int_{\C^{2p}}\frac{e^{-\sum_{k=0}^{2p-1}|\uy_k|^2}}{\pi^{2p}}d\uy
=\frac{1}{p}\sum_{j=1}^{p}\Phi(\zeta_j(\us))=\frac{1}{p}\sum_{j=1}^{p}\Phi(\xi_j)
\end{eqnarray*}
because $(T_{2p}^{(1)}(\us))_j=\xi_j.$

\noindent Let us consider now the case $n>2p$. We cannot use the
same argument used for the case $n=2p$ because $\zeta_j(\us)$ is not
defined for $j=p+1,\dots,n/2$ (see Lemma \ref{lem11}). However by
Lemma \ref{lem1} without loss of generality, we can assume that the
noise is represented by $\mbox{\boldmath $\tilde{\nu}$}_k$ i.e.
$${\bf a}_k=\left\{\begin{array}{ll}\sum_{j=1}^p
c_j\xi_j^k+\sum_{j=p+1}^{n/2}
\tilde{{\bf c}}_{j-p}\tilde{\xi}_{j-p}^k ,\;\; & k=0,\dots,n_p-1 \\
\sum_{j=1}^p c_j\xi_j^k+\mbox{\boldmath $\nu$}_k,\;\; &
k=n_p,\dots,n-1\end{array}\right.$$ where  $n_p=n/2-p.$ We  then
define a new process ${\bf \tilde{a}}_k$ by
$${\bf \tilde{a}}_k=\sum_{j=1}^p c_j\xi_j^k+\mbox{\boldmath $\eta$}_k ,\;\;  k=0,\dots,n-1$$ where
$$\mbox{\boldmath $\eta$}_k=\sum_{j=p+1}^{n/2} \tilde{{\bf c}}_{j-p}\tilde{\xi}_{j-p}^k,$$ and we
consider the process ${\bf a}_k$  as a perturbation of the process
${\bf \tilde{a}}_k$. Let us consider the pencils
$${\bf
P}=[U({\bf a}_1,\dots,{\bf a}_{n-1}),U({\bf a}_0,\dots,{\bf
a}_{n-2})]$$ and
$${\tilde{\bf
P}}=[U({\bf \tilde{a}}_1,\dots,{\bf \tilde{a}}_{n-1}),U({\bf
\tilde{a}}_0, \dots,{\bf \tilde{a}}_{n-2})].$$ We can write
$${\bf P}= {\tilde{\bf P}}+ \sigma{\bf E}$$
where \begin{eqnarray*}{\bf E}=
\frac{1}{\sigma}[U(0,\dots,0,\mbox{\boldmath
$\nu$}_{n_p+1}-\mbox{\boldmath $\eta$}_{n_p+1},
\dots,\mbox{\boldmath $\nu$}_{n-1}-\mbox{\boldmath $\eta$}_{n-1}),
\\ U(0,\dots,0, \mbox{\boldmath
$\nu$}_{n_p}-\mbox{\boldmath $\eta$}_{n_p},\dots,\mbox{\boldmath
$\nu$}_{n-2}-\mbox{\boldmath $\eta$}_{n-2})]=[ {\bf E}_1,{\bf
E}_0].\end{eqnarray*} From \cite{van}, in the limit for
$\sigma\rightarrow 0$, a generalized eigenvalue $\mbox{\boldmath
$\xi$}_j$ of ${\bf P}$ can be expressed as a function of a
generalized eigenvalue $\hat{\xi}_j$ of ${\tilde{\bf P}}$ and
corresponding left and right generalized eigenvectors $\uv_j,\uu_j$
by
\begin{eqnarray*} \mbox{\boldmath $\xi$}_j&=&\hat{\xi}_j+\sigma\frac{\uv_j^H
({\bf E}_1-\hat{\xi}_j{\bf E}_0)\uu_j}{\uv_j^H
{\bf U}_0\uu_j}+O(\sigma^2)\\
&=&\hat{\xi}_j+\sigma\frac{\ue_j^TV^{-1} ({\bf E}_1-\hat{\xi}_j{\bf
E}_0))V^{-T}\ue_j}{{\bf \hat{c}}_j}+O(\sigma^2)\end{eqnarray*}
 where ${\bf U}_0=U({\bf \tilde{a}}_1,\dots,{\bf \tilde{a}}_{n-1})$ and, by
construction,
$$\hat{\xi}_j=\left\{\begin{array}{ll}\xi_j & j=1,\dots,p \\
\tilde{\xi}_{j-p},\;\; & j=p+1,\dots,n/2\end{array}\right . $$
$${\bf \hat{c}}_j=\left\{\begin{array}{ll}c_j & j=1,\dots,p \\
\tilde{{\bf c}}_{j-p},\;\; & j=p+1,\dots,n/2\end{array}\right .
$$
$$V=Vander(\hat{\xi}_1,\dots,\hat{\xi}_{n/2}),\;\;\;
C=diag({\bf \hat{c}}_{1},\dots,{\bf \hat{c}}_{n/2})$$ and
$$\uv_j=\overline{\uu}_j=V^{-H}\ue_j.$$
We notice that we can write $$\ue_j^TV^{-1} ({\bf
E}_1-\hat{\xi}_j{\bf E}_0)V^{-T}\ue_j=\sum_{h=1}^{n/2+p}\gamma_{jh}
{\bf Y}_{n_p+h}$$ where $\gamma_{jh}$ are constants and ${\bf Y}_h($
are i.i.d. zero mean, complex Gaussian variables with unit variance
identified with $\frac{1}{\sqrt{2}\sigma}[\mbox{\boldmath
$\nu$}_{h}-\mbox{\boldmath $\eta$}_{h}],\;h=n_p,\dots,n-1$.

We have \bbs
h_n(z,\sigma)&=&\frac{2}{n}E\left[\sum_{j=1}^{n/2}\delta(z-\mbox{\boldmath
$\xi$}_j)\right]=
\frac{2}{n}E\left[\sum_{j=1}^{p}\delta(z-\mbox{\boldmath
$\xi$}_j)\right]+
\frac{2}{n}E\left[\sum_{j=p+1}^{n/2}\delta(z-\mbox{\boldmath
$\xi$}_j)\right]\\&=& h_n^{(1)}(z,\sigma)+h_n^{(2)}(z,\sigma)\bes By
the same argument used for the case $n=2p$ it follows that
$h_n^{(1)}(z,\sigma)$ converges weakly to
$\frac{2}{n}\sum_{j=1}^p\delta(z-\xi_j)$ when $\sigma\rightarrow 0$.
We then consider $h_n^{(2)}(z,\sigma)$. We have \bbs
h_n^{(2)}(z,\sigma)&=&\frac{2}{n}E\left[\sum_{j=p+1}^{n/2}\delta(z-\mbox{\boldmath
$\xi$}_j)\right]
\\&=&\frac{2}{n}E\left[\sum_{j=p+1}^{n/2}\delta\left(z-\tilde{\xi}_{j-p}-
\sigma\frac{\sum_{h=1}^{n/2+p}\gamma_{jh} {\bf
Y}_{n_p+h}}{\tilde{{\bf c}}_{j-p}}-O(\sigma^2)\right)\right].\bes By
identifying $\frac{\sqrt{n_p}}{\sigma}\tilde{{\bf c}}_{j-p},\;\;
j=p+1,\dots,n/2$ with ${\bf Y}_h,\;\;h=1,\dots,n_p$, which are
i.i.d. zero mean, complex Gaussian variables with unit variance, we
get \bb h_n^{(2)}(z,\sigma)&=& \sum_{j=p+1}^{n/2}\displaystyle
\int_{\C^{n}}\delta\left(z-\tilde{\xi}_{j-p}-
\frac{\sqrt{n_p}}{y_{j-p}}\sum_{h=1}^{n/2+p}\gamma_{jh}
y_{n_p+h}-O(\sigma^2)\right)
\frac{e^{-\frac{1}{\sigma^2}\sum_{k=1}^{n}|y_k|^2}}{\pi^n}d\uy
\nonumber \\&=& \sum_{j=p+1}^{n/2}\displaystyle
\int_{\C^{n-1}}\left[\int_{\C}\delta\left(z-\tilde{\xi}_{j-p}-
\frac{\sqrt{n_p}}{y_{j-p}}\sum_{h=1}^{n/2+p}\gamma_{jh}
y_{n_p+h}-O(\sigma^2)\right)\frac{e^{-|y_{j-p}|^2}}{\pi}dy_{j-p}\right]\nonumber\\&&
\cdot\frac{e^{-\sum_{k=1,k\ne
j-p}^{n}|y_k|^2}}{\pi^{n-1}}d\uy',\;\;\;\{\uy'\}=\{\uy\}\setminus
\{y_{j-p}\}
 \label{eqs}\be
 by making the change of variable $$w=\tilde{\xi}_{j-p}+
\frac{\sqrt{n_p}}{y_{j-p}}\sum_{h=1}^{n/2+p}\gamma_{jh} y_{n_p+h}$$
we get \bbs &&\int_{\C}\delta\left(z-\tilde{\xi}_{j-p}-
\frac{\sqrt{n_p}}{y_{j-p}}\sum_{h=1}^{n/2+p}\gamma_{jh}
y_{n_p+h}-O(\sigma^2)\right)\frac{e^{-|y_{j-p}|^2}}{\pi}dy_{j-p}\\&=&
-\frac{1}{\pi}\int_{\C}\delta\left(z-w-O(\sigma^2)\right)\frac{\sqrt{n_p}
\sum_{h=1}^{n/2+p}\gamma_{jh}
y_{n_p+h}}{(w-\tilde{\xi}_{j-p})^2}e^{-\left|
\frac{\sqrt{n_p}\sum_{h=1}^{n/2+p}\gamma_{jh}
y_{n_p+h}}{w-\tilde{\xi}_{j-p}}\right|^2}dw
\\&=&-\frac{1}{\pi}\frac{\sqrt{n_p}
\sum_{h=1}^{n/2+p}\gamma_{jh}
y_{n_p+h}}{(z-O(\sigma^2)-\tilde{\xi}_{j-p})^2}e^{-\left|
\frac{\sqrt{n_p}\sum_{h=1}^{n/2+p}\gamma_{jh}
y_{n_p+h}}{z-O(\sigma^2)-\tilde{\xi}_{j-p}}\right|^2} .\bes
Inserting this expression in (\ref{eqs}) we get \bbs
h_n^{(2)}(z,\sigma)&=&-\sum_{j=p+1}^{n/2}\frac{\sqrt{n_p}
}{(z-O(\sigma^2)-\tilde{\xi}_{j-p})^2} \displaystyle\\&\cdot&
\sum_{r=1}^{n/2+p}\gamma_{jr}\frac{1}{\pi^{n}}\int_{\C^{n-1}}
y_{n_p+r}e^{-\left| \frac{\sqrt{n_p}\sum_{h=1}^{n/2+p}\gamma_{jh}
y_{n_p+h}}{z-O(\sigma^2)-\tilde{\xi}_{j-p}}\right|^2 -\sum_{k=1,k\ne
j-p}^{n}|y_k|^2}d\uy' \bes and therefore \bbs
\lim_{\sigma\rightarrow
 0}h_n^{(2)}(z,\sigma)&=&-\sum_{j=p+1}^{n/2}\frac{\sqrt{n_p}
}{(z-\tilde{\xi}_{j-p})^2}
\displaystyle\\&\cdot&\sum_{r=1}^{n/2+p}\gamma_{jr}
\frac{1}{\pi^{n}}\int_{\C^{n-1}} y_{n_p+r}e^{-\left|
\frac{\sqrt{n_p}\sum_{h=1}^{n/2+p}\gamma_{jh}
y_{n_p+h}}{z-\tilde{\xi}_{j-p}}\right|^2 -\sum_{k=1,k\ne
j-p}^{n}|y_k|^2}d\uy'=0\bes because \bbs
&&\frac{1}{\pi^{n-1}}\int_{\C^{n-1}} y_{n_p+r}e^{-\left|
\frac{\sqrt{n_p}\sum_{h=1}^{n/2+p}\gamma_{jh}
y_{n_p+h}}{z-\tilde{\xi}_{j-p}}\right|^2 -\sum_{k=1,k\ne
j-p}^{n}|y_k|^2}d\uy'\\&=&\frac{1}{\pi^{n-1}}\int_{\C^{n-1}}y_{n_p+r}e^{-\uy'^H
A \uy'}d\uy'=0,\;\;\mbox{ for a suitable hermitian matrix
}A,\;\;\forall r. \qed \bes

\noindent\underline{Remark}. When the SNR is large the exponential
part dominates the integrand as the Jacobian does not depend on
$\sigma$. Moreover the exponential part has relative maxima close to
$\xi_j$ as expected. In general the integral (\ref{integr}) does not
admit a closed form expression. However when $n=2$, remembering that
the Jacobian with respect to the real and imaginary part of a
complex variable is
 $J_R=|J_C|^2$, the integral (\ref{integr}) becomes
 \begin{eqnarray*} h_2(z,\sigma)&=&\frac{1}{(\pi\sigma^2)^{2}}\int_{\C}\gamma
 e^{-\frac{|\gamma-s_0|^2+|\gamma z-s_1|^2}{\sigma^2}}d\gamma
=\frac{1}{(\pi\sigma^2)^{2}}\int_{\R^2}|\gamma|^2
 e^{-\frac{|\gamma-s_0|^2+|\gamma
 z-s_1|^2}{\sigma^2}}d\Re{\gamma}d\Im{\gamma}\\
 &=&\frac{\sigma^2(1+|z|^2)+|zs_1+s_0|^2}{\pi\sigma^2 (1+|z|^2)^3}
e^{-\frac{|z s_0-s_1|^2}{\sigma^2(1+|z|^2)}}.\end{eqnarray*} We
notice that  $\lim_{\sigma\rightarrow
0}h_2(z,\sigma)=\delta(z-s_1/s_0)=\delta(z-\xi_1)$. Moreover, when
$s_0=s_1=0$ we have $h_2(z,\sigma)=\frac{1}{\pi(1+|z|^2)^2}$ which
is independent of $\sigma^2$, confirming the result obtained in
\cite{barja} for the pure noise case.

\noindent The condensed density has an important role in the
following. Therefore we look for an easily computable approximation.
The following theorem provides a basis for building such an
approximation :
\begin{teo}
 Let be ${\bf F}(z,\overline{z})=
(U_1({\bf \uaa})-zU_0({\bf \uaa}))\overline{(U_1({\bf
\uaa})-zU_0({\bf \uaa}))}$ then
$$E[\log(det\{{\bf F}(z,\overline{z})\})]- \log(det\{E[{\bf F}(z,\overline{z})]\})= o(\sigma)$$
for $\sigma\rightarrow 0$, independently of  z.
 Moreover
\bb E[{\bf
F}(z,\overline{z})]=(U_1(\us)-zU_0(\us))\overline{(U_1(\us)-zU_0(\us))}+
\frac{n\sigma^2}{2}A(z,\overline{z})\label{effetton}\be where
$$A(z,\overline{z})=\left[
\begin{array}{lllll} 1+|z|^2 & \;\;-z &\; 0&\dots&0
\\ -\overline{z}\;\; & 1+|z|^2& \;\;-z &\;\;0&\dots\\ . & . & . & .&.\\
  0 &\dots&0& -\overline{z}\;\; & 1+|z|^2 \end{array}\right].$$
\label{theo3}\end{teo} \noindent\underline{proof}

let us denote by $\mbox{\boldmath $\lambda$}_j$  the eigenvalues of
${\bf F}(z,\overline{z})$ and by $\mu_j$ those of $E[{\bf
F}(z,\overline{z})]$, dropping for simplicity the dependence on
$z,\overline{z}$. Note that $\mu_j\neq E[\mbox{\boldmath
$\lambda$}_j]$, see e.g. \cite[Theorem 8.5]{randpol}. We have
$$E[\log(det\{{\bf F}(z,\overline{z})\})]=\sum_jE[\log(\mbox{\boldmath $\lambda$}_j)]$$ and
$$\log(det\{E[{\bf F}(z,\overline{z})]\})=\sum_j\log(\mu_j)$$
hence it is sufficient to study the difference
$$E[\log(\mbox{\boldmath $\lambda$}_j)]-\log(\mu_j).$$ We then denote by ${\bf \uf}$ the
vector obtained by stacking the real and imaginary parts of the
elements $({\bf F}_{hk},h,k=1,\dots,n/2)$ of ${\bf F}$ and consider
the function
$$g({\bf \uf})=\log(\mbox{\boldmath $\lambda$}_j)$$
and its Taylor expansion around $E[{\bf \uf}]$: \begin{eqnarray*}
g({\bf \uf})&=&g(E[{\bf \uf}])+\sum_{h}\frac{\partial g}{\partial
{\bf \uf}_h}\left|_{{E[{\bf \uf}]}}\right.({\bf \uf}_h-E[{\bf \uf}_h])\\
&+&\frac{1}{2} \sum_{hk}\frac{\partial^2 g}{\partial {\bf
\uf}_h\partial {\bf \uf}_k}\left|_{{E[{\bf \uf}]}}\right.({\bf
\uf}_h-E[{\bf \uf}_h])({\bf \uf}_k-E[{\bf
\uf}_k])+\dots\end{eqnarray*} which can be rewritten as
$$\log(\mbox{\boldmath $\lambda$}_j)-\log(\mu_j)=\sum_{h}\beta_h({\bf \uf}_h-E[{\bf \uf}_h])
+\frac{1}{2} \sum_{hk}\gamma_{hk}({\bf \uf}_h-E[{\bf \uf}_h])({\bf
\uf}_k-E[{\bf \uf}_k])+\dots$$ and, taking expectations,
$$E[\log(\mbox{\boldmath $\lambda$}_j)]-\log(\mu_j)=\frac{1}{2}
\sum_{hk}\gamma_{hk}E[({\bf \uf}_h-E[{\bf \uf}_h])({\bf
\uf}_k-E[{\bf \uf}_k])]+\dots$$
 But
\begin{eqnarray*} {\bf F}(z,\overline{z})
&=&(U_1(\us)-zU_0(\us)\overline{(U_1(\us)-zU_0(\us))}\\
&+&(U_1(\mbox{\boldmath $\unu$})-zU_0(\mbox{\boldmath $\unu$}))\overline{(U_1(\mbox{\boldmath $\unu$})-zU_0(\mbox{\boldmath $\unu$}))}\\
&-&(U_1(\us)-zU_0(\us)\overline{(U_1(\mbox{\boldmath $\unu$})-zU_0(\mbox{\boldmath $\unu$}))}\\
&-&(U_1(\mbox{\boldmath $\unu$})-zU_0(\mbox{\boldmath
$\unu$}))\overline{(U_1(\us)-zU_0(\us))}
\end{eqnarray*}
and
\begin{eqnarray*} E[{\bf F}(z,\overline{z})]
&=&(U_1(\us)-zU_0(\us)\overline{(U_1(\us)-zU_0(\us))}\\
&+&E[(U_1(\mbox{\boldmath $\unu$})-zU_0(\mbox{\boldmath $\unu$}))\overline{(U_1(\mbox{\boldmath $\unu$})-zU_0(\mbox{\boldmath $\unu$}))}]\\
&=&(U_1(\us)-zU_0(\us)\overline{(U_1(\us)-zU_0(\us))}+
\frac{n\sigma^2}{2}A(z,\overline{z})
\end{eqnarray*}
by a straightforward computation similar to that given in
\cite[Th.3]{barja} for the pure noise case. Therefore
\begin{eqnarray*}{\bf F}(z,\overline{z})-E[{\bf F}(z,\overline{z})]&=&
(U_1(\mbox{\boldmath $\unu$})-zU_0(\mbox{\boldmath $\unu$}))\overline{(U_1(\mbox{\boldmath $\unu$})-zU_0(\mbox{\boldmath $\unu$}))}\\
&-&(U_1(\us)-zU_0(\us)\overline{(U_1(\mbox{\boldmath $\unu$})-zU_0(\mbox{\boldmath $\unu$}))}\\
&-&(U_1(\mbox{\boldmath $\unu$})-zU_0(\mbox{\boldmath $\unu$}))\overline{(U_1(\us)-zU_0(\us))}\\
&-&\frac{n\sigma^2}{2}A(z,\overline{z})
\end{eqnarray*}
hence $E[({\bf \uf}_h-E[{\bf \uf}_h])({\bf \uf}_k-E[{\bf \uf}_k])]$
is a linear combination of functions of $z$ and $\overline{z}$ with
coefficients equal to either $\sigma^2$ or $\sigma^4$ because the
odd moments of a Gaussian are zero. By a similar argument all the
dropped terms in the Taylor expansion above will depend on even
powers of $\sigma$. Hence
$$E[\log(\mbox{\boldmath $\lambda$}_j)]-\log(\mu_j)=o(\sigma)$$ independently of
$z,\overline{z}.\qed$

By noticing that $|{\bf Q}(z)|^2=det\{{\bf F}(z,\overline{z})\}$, an
approximation of the condensed density is then given by
$$\tilde{h}_n(z,\sigma)=\frac{1}{2\pi n}\Delta\sum_{\mu_j(z)> 0
}\log(\mu_j(z))$$ where $\mu_j(z)$ are the eigenvalues of $E[{\bf
F}(z,\overline{z})].$  Unfortunately $\tilde{h}_n(z,\sigma)$ is not
a probability density as it can eventually  assume negative values.
However the following results hold

\begin{teo}
The function $\tilde{h}_n(z,\sigma)$ is continuous in $\sigma$ and
in $z$. In the limit cases $\sigma=0$ and $\{c_k=0,k=1,\dots,p\}$ it
is  given respectively by
$$\tilde{h}_n(z,0)=\frac{2}{n}\sum_{j=1}^p\delta(z-\xi_j)$$ and by
$$\tilde{h}_n(z,\sigma)= \frac{1}{4\pi}\Delta w_n(z)$$
where $$w_n(z)=\frac{1}{n}\log\sum_{j=0}^n|z|^{2j}.$$ Moreover, in
this second case,
$\lim_{n\rightarrow\infty}\tilde{h}_n(z,\sigma)=\delta(|z|-1).$\label{teo4}
\end{teo}
\noindent\underline{proof}

 $\tilde{h}_n(z,\sigma)$ is  continuous in $\sigma$
and in $z$ because of the continuous dependence of the eigenvalues
on the elements of the corresponding  matrix. When $\sigma=0$, let
$V\in\C^{n/2,p}$ be the Vandermonde matrix such that
$U_0(\us)=VCV^T$ and $U_1(\us)=VCZV^T$. Let $V=QR$ be the $QR$
decomposition of  $V$. Then
$$E[{\bf F}(z,\overline{z})]=QRC(Z-zI)R^TQ^T\overline{Q}
\overline{R}\overline{(Z-zI)}CR^HQ^H.$$ But
$R=\left(\begin{array}{c}\tilde{R} \\ 0 \end{array}\right),$
therefore $R^T\overline{R}=\tilde{R}^T\overline{\tilde{R}}$;
moreover $Q^T\overline{Q}=I$, hence the eigenvalues of $E[{\bf
F}(z,\overline{z})]$ are the same of those of the matrix
$$RC(Z-zI)R^T\overline{R}\overline{(Z-zI)}CR^H=
\left(\begin{array}{cc}\tilde{R}C(Z-zI)\tilde{R}^T\overline{\tilde{R}}
\overline{(Z-zI)}C\tilde{R}^H & \;\;\;\;0 \\ \;\;\;\;0 &\;\; \;\;0
\end{array}\right).$$  The non-zero
eigenvalues of $E[{\bf F}(z,\overline{z})]$ are then the same of
those of the matrix
$$\tilde{R}C(Z-zI)\tilde{R}^T\overline{\tilde{R}}
\overline{(Z-zI)}C\tilde{R}^H.$$ We then have
\begin{eqnarray*} \tilde{h}_n(z,0)
&=&\frac{1}{2\pi n}\Delta\sum_{\mu_j(z)> 0}\log(\mu_j(z))\\
&=&\frac{1}{2\pi n}\Delta\log\left(\prod_{j=1}^p|z-\xi_j|^2\cdot
|det(\tilde{R})|^4\prod_{j=1}^p c_j^2\right)\\
&=&\frac{2}{4\pi n}\sum_{j=1}^p\Delta\log|z-\xi_j|^2
=\frac{2}{n}\sum_{j=1}^p\delta(z-\xi_j)
\end{eqnarray*} because
$\frac{1}{4\pi}\Delta\log(|z|^2)=\delta(z)$ (see e.g.
\cite[pg.47]{sch}).
 \noindent When
$\{c_k=0,k=1,\dots,p\}$
\begin{eqnarray*}
\tilde{h}_n(z,\sigma)&=&\frac{1}{2\pi
n}\Delta\log(det\{A(z,\overline{z})\}) =\frac{1}{2\pi
n}\Delta\log(\sum_{j=0}^n|z|^{2j}).
\end{eqnarray*}
The last part of the thesis follows by the same argument used in the
proof of  Theorem 3 in \cite{barja}. $\qed$

\begin{cor}
$\tilde{h}_n(z,\sigma)-h_n(z,\sigma)$ converges weakly to $0$ when
$\sigma\rightarrow 0$
\end{cor}
\noindent\underline{proof}

Let $\Phi(z)$ be a nonnegative test function supported on ${\C}$.
Denoting by $h_n^*(z)=\frac{2}{n}\sum_{j=1}^p\delta(z-\xi_j)$, from
Theorems \ref{teo2} and  \ref{teo4} we have $\forall
\nu>0,\;\;\exists \sigma_1  \mbox{ and } \sigma_2>0$ such that
$$\left|\int_{\C} \Phi(z)\left(h_n(z,\sigma)-h_n^*(z)
\right)dz\right|<\frac{\nu}{2},\;\;\forall \sigma<\sigma_1$$ and
$$\left|\int_{\C} \Phi(z)\left(\tilde{h}_n(z,\sigma)-h_n^*(z)
\right)dz\right|<\frac{\nu}{2},\;\;\forall \sigma<\sigma_2$$ hence,
if $\sigma_\nu=\min\{\sigma_1,\sigma_2\}$, we have $\forall
\sigma<\sigma_\nu$
\begin{eqnarray*}& &\left|\int_{\C} \Phi(z)\left(h_n(z,\sigma)-
\tilde{h}_n(z,\sigma)\right)dz\right| \\ &\le& \left|\int_{\C}
\Phi(z)\left(h_n(z,\sigma)-h^*(z) \right)dz\right|+ \left|\int_{\C}
\Phi(z)\left(\tilde{h}_n(z,\sigma)-h^*(z)
\right)dz\right|\le\nu.\qed\end{eqnarray*}

\section{Identifiability of $S(z)$ and approximation properties of
$E[{\bf S}_n(z)]$}

We want now to exploit the information about the location in the
complex plane of the Pade' poles, provided by the condensed density
$h_n(z)$, to prove some properties relating ${\bf
S}_n(z)=\sum_{j=1}^{n/2}{\bf c}_j\delta(z-\mbox{\boldmath $\xi$}_j)$
to the true measure $S(z)$.

Before affording the problem of estimating $S(z)$ from the data
${\bf \uaa}$ we need to check that the data provide enough
information to solve it. Precise conditions that must be met to
solve the problem are well known in the noiseless case and are
reported in the introduction. When noise is present the
identifiability problem is an open one. Its solvability can depend
on the amount of "a priori" information available \cite{donoho}
and/or on the ability to devise smart algorithms. In the following a
definition of identifiability is given and, based on it, some
properties of ${\bf S}_n(z)$ are proved.

\begin{defi}
The measure $S(z)$  is identifiable from the data ${\bf
a_k},k=0,\dots,n-1$ if $\exists \;\;r_k>0,k=1,\dots,p$ such that
\begin{itemize}
\item
$h_n(z) \mbox{ is unimodal in }  N_k=\{z|\;|z-\xi_k|\le r_k\}$
\item
$ \bigcap_{k=1}^{p}N_k=\emptyset$
\end{itemize}

\end{defi}
The idea is that $S(z)$ can be identified from the data ${\bf \uaa}$
if the random generalized eigenvalues have a condensed density with
separate peaks centered on $\xi_j,j=1,\dots,p$. As, by Theorem
\ref{teo2}, $h_n(z,\sigma)$ converges weakly to
$\frac{2}{n}\sum_{j=1}^p\delta(z-\xi_j)$ when $\sigma\rightarrow 0$,
it must exists a $\sigma'>0$ small enough to make $S(z)$
identifiable $\forall \; \sigma<\sigma'$.

In order to apply the proposed method one should check that the
identifiability conditions are verified. As $h_n(z,\sigma)$ depends
on the unknown quantities $p,c_j,\xi_j$ this is of course
impossible. However in most real problems we have some prior
information about the unknown measure $S(z)$ that we can exploit to
get reasonable interval estimates for $p,c_j,\xi_j$. Moreover in
many instances either $n$ or $\sigma$ or both can be freely chosen.
By Theorem \ref{theo3}, equation \ref{effetton},  $n$ should not be
as large as possible to get the best estimates of $S(z)$. In fact
too many data will convey too much noise which could mask the signal
$s_k$. We can therefore properly design an experiment by computing
$h_n(z,\sigma)$ for many values of $n$ and $\sigma$ and choose
$n_{ott}$ and $\sigma_{ott}$ (optimal design) that make identifiable
the  measures corresponding to prior estimates of $p,c_j,\xi_j$. To
identify the unknown measure $S(z)$ we then hopefully need to
measure $n_{ott}$ data affected by an error with s.d.
$\sigma_{ott}$. Unfortunately $h_n(z)$ does not admit a closed form
expression and to compute the expectation that appears in its
definition  we need to perform a time consuming MonteCarlo
experiment. This is why we proposed an approximation
$\tilde{h}_n(z)$ of $h_n(z)$ which can be quickly computed by
solving hermitian eigenvalues problems.

Let us consider the function
$$S_n(z)=E[{\bf
S}_n(z)]=\sum_{j=1}^{n/2}E[{\bf c}_j\delta(z-\mbox{\boldmath
$\xi$}_j]$$ where $\{{\bf c}_j,\mbox{\boldmath
$\xi$}_j\},\;j=1,\dots,n/2\}$ are the solution of the complex
exponential interpolation problem for the data $\{{\bf
a_k},k=0,\dots,n-1\}$.

The relation between $S_n(z)$ and the unknown measure $S(z)$ is
given by the following
\begin{teo}
If $S(z)$ is identifiable from ${\bf \uaa}$ then
$$\int_{N_h}S_n(z)dz= c_h+o(\sigma)$$
and $$\int_{A}S_n(z)dz= o(\sigma),\;\;\forall A\subset
D-\bigcup_jN_j$$\label{teo33}
\end{teo}

\noindent\underline{proof}

From the identifiability hypothesis we know that
$$\int_{N_k}h_n(z)dz=\frac{2}{n}\sum_{j=1}^{n/2}Prob[\mbox{\boldmath $\xi$}_j\in
N_k]>0,\;\;k=1,\dots,p.$$ Therefore there exist $\mbox{\boldmath
$\xi$}_{j_k}$ such that $Prob[\mbox{\boldmath $\xi$}_{j_k}\in
N_k]>0$. Among the $\mbox{\boldmath $\xi$}_{j_k}$ let us denote by
$\mbox{\boldmath $\xi$}_{\hat{k}}$ the one such that
$Prob[\mbox{\boldmath $\xi$}_{j_k}\in N_k]$ is maximum. From the
identifiability hypothesis the $\mbox{\boldmath $\xi$}_{\hat{k}}$
are distinct. Moreover all the
$\mbox{\boldmath$\xi$}_j,\;j=1,\dots,n/2 $ can be sorted in such a
way that
$\mbox{\boldmath$\xi$}_j=\mbox{\boldmath$\xi$}_{\hat{j}},\;j=1,\dots,p$
and, by Lemma \ref{lem11}, to $\mbox{\boldmath$\xi$}_k$ it
corresponds ${\bf c}_k$ such that
$$E[{\bf c}_k]=\left\{\begin{array}{ll}c_k +o(\sigma),& k=1,\dots,p \\
o(\sigma),\;\; & k=p+1,\dots,n/2\end{array}\right . $$ But then for
$k=1,\dots,p$
$$\int_{N_k}S_n(z)dz=
\sum_{j=1}^{n/2}\int_{N_k}E[{\bf c}_j\delta(z-\mbox{\boldmath
$\xi$}_j]dz=$$
$$=\sum_{j=1}^{n/2}\int_{N_k}\left(\int_{
\C^2}\gamma\delta(z-\zeta)d\mu_{\gamma \zeta}\right)dz=$$
$$=\sum_{j=1}^{n/2}\int_{
\C^2}\gamma\left(\int_{N_k}\delta(z-\zeta)dz\right)d\mu_{\gamma
\zeta}$$ where $\mu_{\gamma \zeta}$ is the joint distribution of
${\bf c}_j$ and $\mbox{\boldmath $\xi$}_j$. We have
$$\int_{N_k}\delta(z-\zeta)dz=\left\{\begin{array}{ll}
1 & \mbox{if $\zeta\in N_k$} \\
0 & \mbox{otherwise}\end{array}\right.$$ hence,
$$\int_{N_k}S_n(z)dz=\sum_{j=1}^{n/2}
E[{\bf c}_j \delta_{jk}]=E[{\bf c}_k]=c_k+o(\sigma) .$$ By a similar
argument the second part of the thesis follows. $\qed$

\section{The $\mathsf{P}$-transform}

In order to solve the original moment problem we need to compute
$$S_n(z,\sigma^2)=\sum_{j=1}^{n/2}E[{\bf c}_j\delta(z-\mbox{\boldmath $\xi$}_j].$$ In
order to estimate the expected value we build independent
replications of the data (pseudosamples) by defining
$${\bf a}_k^{(r)}={\bf a}_k+\mbox{\boldmath $\nu$}_k^{(r)},\;\;k=0,\dots,n-1;\;\;\;r=1,\dots,R$$
where $\{\mbox{\boldmath $\nu$}_k^{(r)}\}$ are i.i.d. zero mean
complex Gaussian variables with variance $\sigma'^2$ independent of
${\bf a}_h,\;\forall h$. Therefore
$$E[{\bf a}_k^{(r)}]=s_k,\;\;\;E[({\bf a}_k^{(r)}-s_k)
({\bf \overline{a}}_h^{(s)}-\overline{s}_h)]=
\tilde{\sigma}^2\delta_{hk}\delta_{rs}$$ where
$\tilde{\sigma}^2=\sigma^2+\sigma'^2$.
 For $r=1,\dots,R$, we
define the statistics
$${\bf \hat{S}}_{n,r}(z,\tilde{\sigma}^2)=\sum_{j=1}^{n/2}
{\bf c}_j^{(r)}\delta(z-\mbox{\boldmath $\xi$}_j^{(r)})   $$ where
${\bf c}_j^{(r)}, \mbox{\boldmath $\xi$}_j^{(r)}$ are the solution
of the complex exponentials interpolation problem for the data ${\bf
a}_k^{(r)},\;\;k=0,\dots,n-1.$ As, by Lemma \ref{lem11}, the
transformation
$$T:\{{\bf a}_k^{(r)},k=0,\dots,n-1\}\rightarrow \{[{\bf c}_j^{(r)},
\mbox{\boldmath $\xi$}_j^{(r)}],j=1,\dots,n/2\}$$ is one-to-one,
${\bf \hat{S}}_{n,r}(z,\tilde{\sigma}^2)$ are i.i.d. with mean
$S_n(z,\tilde{\sigma}^2)$ and finite variance
$\zeta(z,\tilde{\sigma}^2)$ because $\{\mbox{\boldmath
$\nu$}_k^{(r)}\}$ are i.i.d. . Therefore the statistic
$${\bf \hat{S}}_{n,R}(z,\tilde{\sigma}^2)=\frac{1}{R}\sum_{r=1}^R
{\bf \hat{S}}_{n,r}(z,\tilde{\sigma}^2)$$ has mean
$S_n(z,\tilde{\sigma}^2)=E[{\bf \hat{S}}_{n,r}(z,\tilde{\sigma}^2)]$
and variance $\frac{1}{R}\zeta(z,\tilde{\sigma}^2).$

\noindent Let us consider the statistic
$${\bf \hat{S}}_{n}(z,\sigma^2)=
\sum_{j=1}^{n/2}{\bf c}_j\delta(z-\mbox{\boldmath $\xi$}_j) $$ where
${\bf c}_j, \mbox{\boldmath $\xi$}_j$ are the solution of the
complex exponentials interpolation problem for the data ${\bf
a}_k,\;\;k=0,\dots,n-1$ and the conditioned statistic
$${\bf \hat{S}}_{n,R}^c(z,\tilde{\sigma}^2)={\bf \hat{S}}_{n,R}(z,\tilde{\sigma}^2)|{\bf \uaa}$$ which
are both computable from the observed data $\uaa$. We have
\begin{lemma} For $n$ and $\sigma>0$ fixed and $\forall z$ and $\tilde{\sigma}$,
 $$E[{\bf \hat{S}}_{n,R}^c(z,\tilde{\sigma}^2)]=S_n(z,\tilde{\sigma}^2)$$
$$\lim_{R\rightarrow\infty}var[{\bf \hat{S}}_{n,R}^c(z,\tilde{\sigma}^2)]=0.$$ \label{lemma4}
\end{lemma}
\noindent\underline{proof}

\noindent from the conditional variance formula (\cite{reny}) we
have
$$E[{\bf \hat{S}}_{n,R}^c(z,\tilde{\sigma}^2)]=E[{\bf \hat{S}}_{n,R}(z,\tilde{\sigma}^2)]=
S_n(z,\tilde{\sigma}^2)$$ and
$$var[({\bf \hat{S}}_{n,R}^c(z,\tilde{\sigma}^2)]\le
var[{\bf
\hat{S}}_{n,R}(z,\tilde{\sigma}^2)]=\frac{1}{R}\zeta(z,\tilde{\sigma}^2).\;\;\qed$$

\noindent It follows that $\forall z$ the risk of ${\bf
\hat{S}}_{n,R}^c(z,\tilde{\sigma}^2)$ as an estimator of $S(z)$ with
respect to the loss function given by the absolute difference could
be smaller than the risk of the estimator ${\bf
\hat{S}}_{n}(z,\sigma^2)$ if $R$ and $\tilde{\sigma}$ are suitably
chosen, despite of the fact that its bias is larger because
$\tilde{\sigma}>\sigma$ and Theorem \ref{teo33} holds. As a matter
of fact this possibility is always verified provided that $\sigma'$
and $R$ are suitably chosen as proved in the following
\begin{teo}
Let $M(z)$ and $M_c(z)$ be the mean squared error of ${\bf
\hat{S}}_{n}(z,\sigma^2)$ and ${\bf
\hat{S}}_{n,R}^c(z,\tilde{\sigma}^2)$ respectively. In the limit for
$\sigma\rightarrow 0$, it exist $\sigma'$ and $R(\sigma')$ such that
$\forall R\ge R(\sigma')$, $M_c(z)<M(z)\;\; \forall z$.
\end{teo}
\noindent\underline{proof}

let $M_c(z)=v_c+b_c^2$ and $M(z)=v+b^2$ be the decomposition of the
mean squared errors in the sum of variance plus squared bias. Then
$M_c(z)-b^2=v_c+(b_c^2-b^2)$. By Lemma \ref{lemma4}, $b_c$  is equal
to the bias of ${\bf \hat{S}}_{n}(z,\tilde{\sigma}^2)$ and, by
Theorem \ref{teo33}, it is  $o(\tilde{\sigma})$  for
$\tilde{\sigma}\rightarrow 0$. Then $\lim_{\sigma'\rightarrow
0^+}(b_c^2-b^2)=0$. Moreover, by Lemma \ref{lemma4},
$\lim_{R\rightarrow\infty}v_c=0$. Therefore $\forall v>0,\;\;
\exists \sigma'_v$ and $R(\sigma'_v)$ such that $\forall
\sigma'<\sigma_v',\;\;v_c+(b_c^2-b^2)<v$ and then
$M_c(z)<M(z).\;\;\qed$

\noindent In order to define a discrete transform,
  we evaluate ${\bf \hat{S}}_{n,R}^c(z,\tilde{\sigma}^2))$ on a
lattice $L=\{(x_i,y_i),i=1,\dots,N\}$ such that
$$\min_j\Re{\xi_j}>\min_i x_i;\;\;\;\max_j\Re{\xi_j}<\max_i x_i$$
$$\min_j\Im{\xi_j}>\min_i y_i;\;\;\;\max_j\Im{\xi_j}<\max_i y_i.$$
In order to cope with the Dirac distribution appearing in the
definition of ${\bf \hat{S}}_{n,R}^c(z,\tilde{\sigma}^2))$ it is
convenient to use an alternative expression given by
$${\bf \hat{S}}_{n,R}^c(z,\tilde{\sigma}^2)=\frac{1}{2\pi R}\Delta\left(\sum_{r=1}^R\sum_{j=1}^{n/2}
[{\bf c}_j^{(r)}|{\bf \uaa}]\log(|z-[\mbox{\boldmath
$\xi$}_j^{(r)}|{\bf \uaa}]|) \right)$$ which can be obtained by the
former one by remembering that
$\frac{1}{4\pi}\Delta\log(|z|^2)=\delta(z)$ (see e.g.
\cite[pg.47]{sch}). In this way the problem of discretizing the
Dirac $\delta$ is reduced to  discretizing the Laplacian operator,
which is easier to cope with. We then get  a random matrix
$\mathsf{P}(\tilde{\sigma}^2)\in\Re_+^{(N\times N)}$ such that
$\mathsf{P}(h,k,\tilde{\sigma}^2)={\bf \hat{S}}_{n,R}^c(x_h+iy_k)$.
We call this matrix the $\mathsf{P}$-transform of the vector $[{\bf
a}_0,\dots,{\bf a}_{n-1}]$.

\section{Estimation procedure}

The $\mathsf{P}$-transform gives a global picture of the measure
$S(z)$. However an estimate of the unknown parameters $p$, $\{\xi_j,
c_j,j=1,\dots,p\}$ are usually of interest. An automatic procedure
to get such estimates is now described. Let $
\mathsf{P}(\tilde{\sigma}^2)$ be the $\mathsf{P}$-transform computed
by using  $R$ pseudosamples with variance $\tilde{\sigma}^2$. The
proposed procedure is the following (dropping for simplicity the
conditioning to ${\bf \uaa}$):
\begin{itemize}
\item memorize all the Pade' poles $\mbox{\boldmath $\xi$}_j^{(r)}$
 and the corresponding
residuals ${\bf c}_j^{(r)},\;\;r=1,\dots,R$ used for computing
$\mathsf{P}(\tilde{\sigma}^2)$
\item identify the local maxima of
$\mathsf{P}(\tilde{\sigma}^2)$  and sort them in increasing order
with respect to the local maxima values. The local maxima are
candidate estimates of $\{\xi_j,j=1,\dots,p\}$
\item for each candidate a cluster of (previously
memorized) Pade' poles was estimated by including all the poles
closest to the current candidate until the cluster cardinality
equals a predefined percentage (e.g.$>50\%$) of the number $R$ of
pseudosamples. The rationale is that if the candidate is close to
one of the  $\xi_j$ most of the pseudosamples should provide a Pade'
pole close to it. Notice that spurious  clusters - i.e. not centered
close to some $\xi_j\;$- can be expected \cite{bama1}
\item
all the candidates whose associated cluster does not have the
prescribed cardinality are eliminated.
 The number
$\hat{p}$ of left candidates is then an estimate of $p$
\item
for each of the $\hat{p}$ clusters the Pade' poles and the
corresponding residuals (previously memorized) were then averaged
and provided estimates $\hat{\xi}_j,\hat{c}_j,j=1,\dots,\hat{p}$ of
the unknown parameters. Hopefully to $\hat{\xi}_j$ associated to
spurious clusters should correspond relatively small $\hat{c}_j$.
\end{itemize}

\section{Numerical results}
\begin{figure}[h]
\begin{center}
\hspace{1.7cm}\centerline{\fbox{\epsfig{file=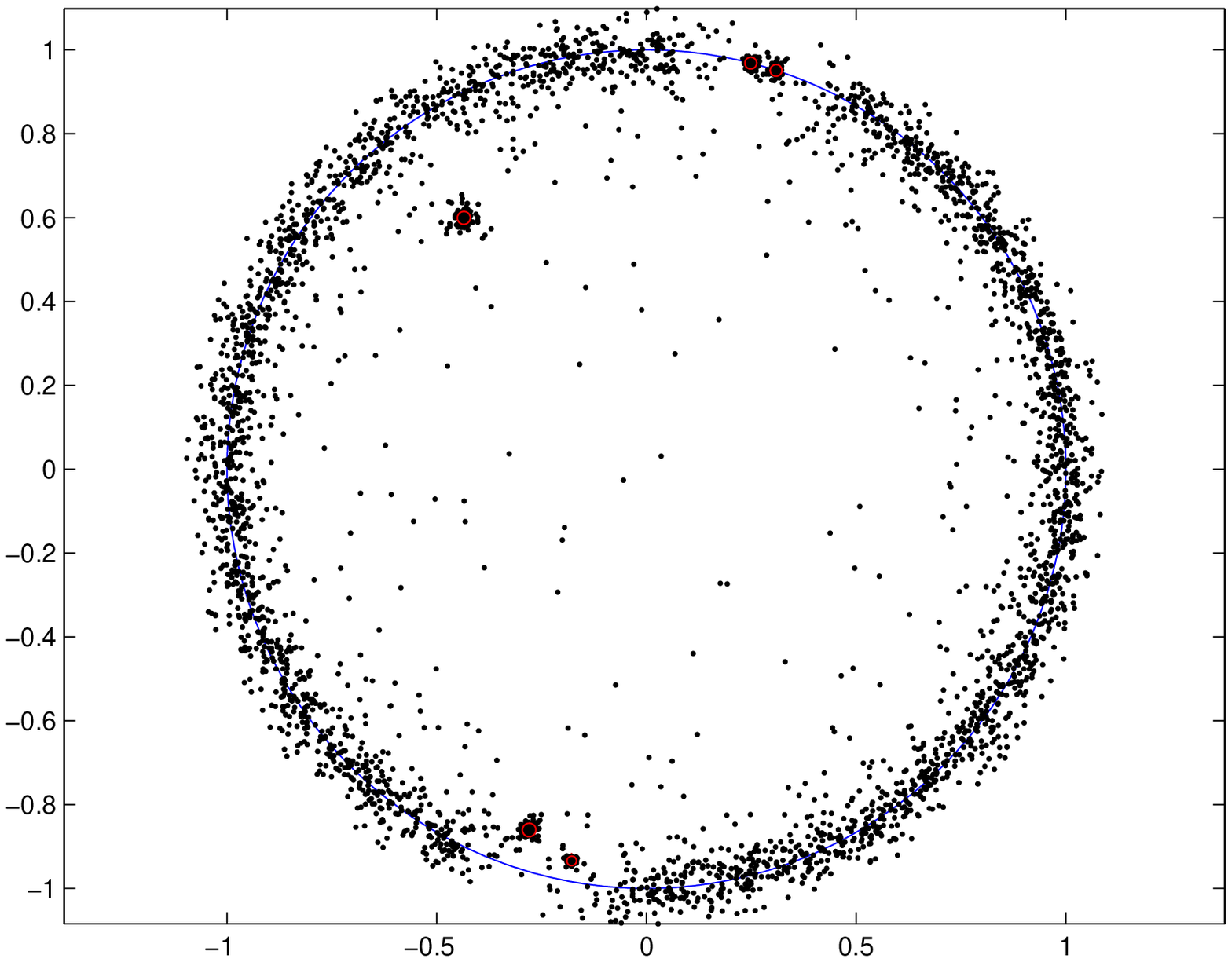,height=4.8cm,width=6cm}}
\hspace{.5cm}\fbox{\epsfig{file=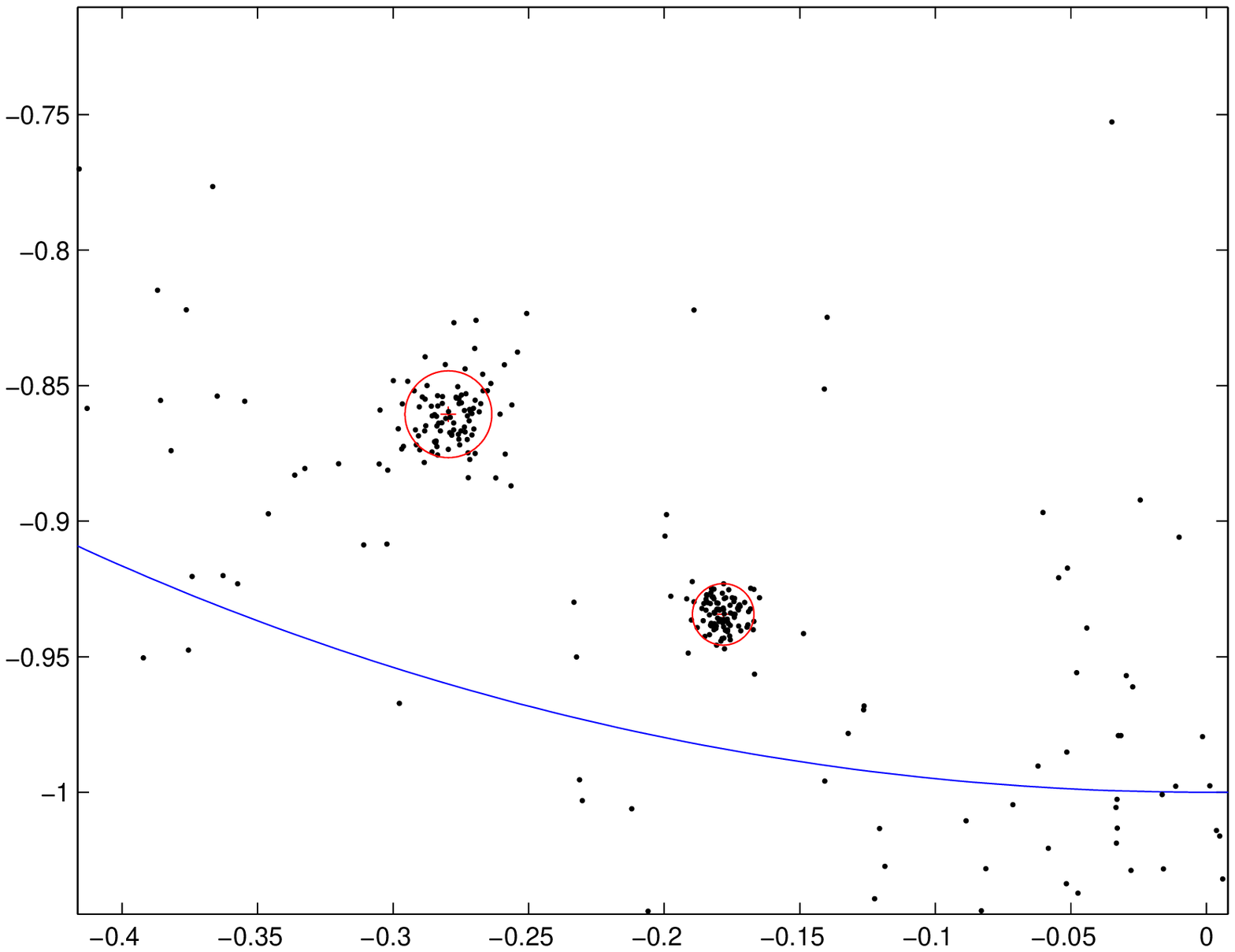,height=4.8cm,width=6cm}}}
\vspace{.2cm}
\hspace{1.7cm}\centerline{\fbox{\epsfig{file=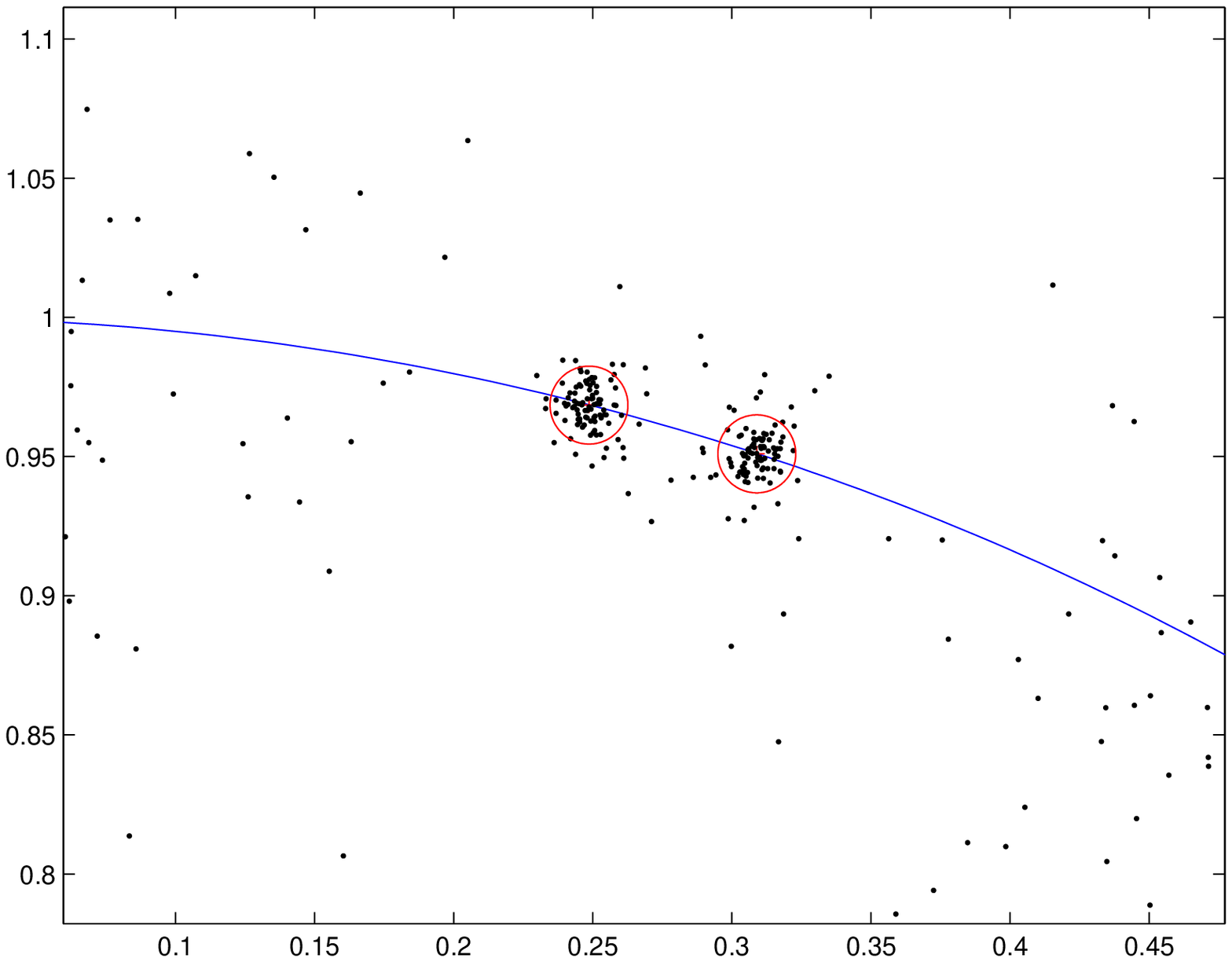,height=4.8cm,width=6cm}}
\hspace{.5cm}\fbox{\epsfig{file=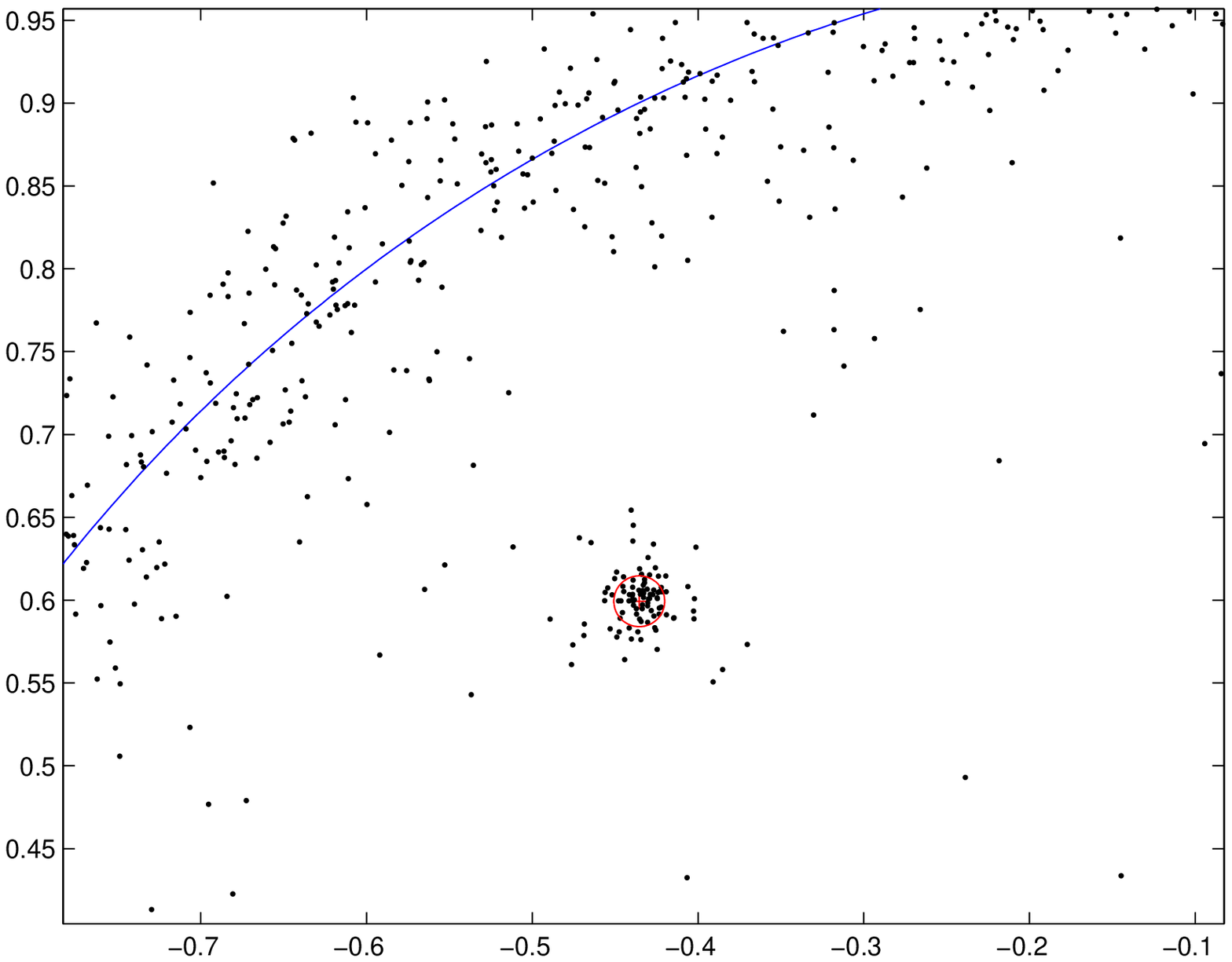,height=4.8cm,width=6cm}}}
\end{center}
\caption{Top left: location of Pade' poles for $100$ independent
realizations of the noise; the circles are the estimated support of
the condensed density in a neighborhood of $\xi_j$; top right:zoom
in a neighborhood of the 1-st and 2-nd components; bottom left: zoom
in a neighborhood of the 3-rd and 4-th components; zoom in a
neighborhood of the 5-th component (see section 4).}  \label{fig1}
\end{figure}

In this section some experimental evidence of the claims made in the
previous sections is given.  A model with $p=5$ components given by
$$\underline{\xi}=\left[ e^{-0.1-i 2\pi  0.3},e^{-0.05-i 2\pi
0.28},e^{-0.0001+i 2\pi 0.2},e^{-0.0001+i 2\pi  0.21},e^{-0.3-i 2\pi
0.35}\right]$$ $$ \underline{c}=\left[ 6,3,1,1,20\right
],\;\;\sigma=0.2,\;\;n=80$$ is considered. We notice that $SNR=5$
and the frequencies of the $3^{rd}$ and $4^{th}$ components are
closer than the Nyquist frequency ($0.21-0.20=0.01<1/n=0.0125$).
Hence a superesolution problem is involved in this case. The quality
of the approximation of $\tilde{h}(z)$ to the condensed density is
first addressed, $\tilde{h}(z)$ is then computed along a line which
pass through $\xi_j$ and the closest among the $(\xi_h,h\ne j)$. If
the model is identifiable
 $\tilde{h}(z)$ should have a local maximum close to
$\xi_j$ along this line. The interquartile range $\hat{r}_j$ of a
restriction of $\tilde{h}(z)$ to a neighbor of this maximum is then
considered as an estimate of the radius of the local support of
$\tilde{h}(z)$ assumed circular. Then $M=100$ independent data sets
$\uaa^{(m)}$ of length $n$  were generated and the Pade' poles
$\underline{\xi}^{(m)},m=1,\dots,M$ were plotted in fig.\ref{fig1}
where  circles of radii $\hat{r}_j$ centered on $\xi_j$ have been
represented too. We notice that the
 circles are reasonable estimates of the Pade' poles clusters
which provide an estimate of the support of the peaks of the true
condensed density corresponding to $\xi_j,j=1,\dots,p$. We conclude
that $\tilde{h}(z)$ is a reliable approximation of the condensed
density and therefore, with the choice of $n$ and $\sigma$ made
above, the model is likely to be identifiable.

We want now to show by means of a small simulation study the quality
of the estimates of the parameters $\underline{\xi}$ and $
\underline{c}$ which define the unknown measure $S(z)$. To this aim
the bias, variance and mean squared error (MSE) of each parameter
separately will be estimated.  $M=500$ independent data sets
$\uaa^{(m)}$ of length $n$  were generated by using the model
parameters given above. For $m=1,\dots,M$  the
$\mathsf{P}$-transform $\mathsf{P}^{(m)}$ was computed based on
$R=100$ pseudosamples with $\sigma'^2=10^{-4}\sigma^2$ on a square
grid of dimension $N=200$. The estimation procedure is then applied
to each of the $\mathsf{P}^{(m)},m=1,\dots,M$ and the corresponding
estimates
$\hat{\xi}_j^{(m)},\hat{c}_j^{(m)},j=1,\dots,\hat{p}^{(m)}$ of the
unknown parameters were obtained. As we know the true value $p$, if
less than $p$ local maxima were found in the second step or if
$\hat{p}^{(m)}<p$ in the fourth step of the procedure, the
corresponding data set $\uaa^{(m)}$ was discarded.

In Table 1 the bias, variance and MSE of each parameter including
$p$ is reported. They were computed by choosing among the
$\hat{\xi}_j^{(m)},j=1,\dots,\hat{p}^{(m)}$ the one closest to each
$\xi_k,k=1,\dots,p$ and the corresponding $\hat{c}_j^{(m)}$. If more
than one $\xi_k$ is estimated by the same $\hat{\xi}_j^{(m)}$ the
$m-$th data set $\uaa^{(m)}$ was discarded. In the case considered
$65\%$ data sets were accepted. Looking at Table 1 we can conclude
that the true measure can be estimated quite accurately  in $65\%$
of cases.

When $\hat{p}_j^{(m)}>p$ we computed also the average residual
amplitude
$$a_{res}=\frac{1}{|\tilde{M}|}\sum_{m\in\tilde{M}}
\frac{1}{(\hat{p}^{(m)}-p)}\sum_{j=p+1}^{\hat{p}^{(m)}}\hat{c}_j^{(m)},\mbox{
where } \tilde{M}=\{m|\hat{p}_j^{(m)}>p\}$$ which represents the
contribution to ${\bf \hat{S}}_{n,R}^c(z,\tilde{\sigma}^2))$ of all
the components which give rise to spurious clusters. In the case
considered its value is $a_{res}= 1.165$ which should be compared
with the true amplitudes $\underline{c}$. We can conclude that even
when more components then the true ones are detected their relative
importance is very low.

In order to appreciate the  advantage of the estimator ${\bf
\hat{S}}_{n,R}^c(z,\tilde{\sigma}^2)$ with respect to ${\bf
\hat{S}}_{n}(z,\sigma^2)$, the same $M=100$ independent data sets
$\uaa^{(m)}$ of length $n$ generated before were considered. The
corresponding Pade' poles and weights
$(\hat{\xi}_j^{(m)},\hat{c}_j^{(m)},j=1,\dots,n/2)$ were computed
and ordered for each $m$ in decreasing order w.r. to the absolute
value of the weights. The true $(\xi_j,c_j,j=1,\dots,p)$ were
ordered in the same way and the error
$$e_0(m)=\sum_{j=1}^p(\hat{\xi}_j^{(m)}-\xi_j)^2+\sum_{j=1}^p(\hat{c}_j^{(m)}-c_j)^2  $$
was computed for $m=1,\dots,M$ and plotted in fig.2. Then to each of
the $M$ data sets $\uaa^{(m)}$ previously generated $R=100$ i.i.d.
zero-mean Gaussian samples with variance $\sigma'^2=0.64\sigma^2$
were added and
$(\hat{\xi}_j^{(m,r)},\hat{c}_j^{(m,r)},j=1,\dots,n/2,\;\;r=1,\dots,R)$
were computed and ordered as before for each $m$ and $r$. Finally
the error
$$e_R(m)=\sum_{j=1}^p\left(\frac{1}{R}\sum_{r=1}^R\hat{\xi}_j^{(m,r)}-\xi_j\right)^2+
\sum_{j=1}^p\left(\frac{1}{R}\sum_{r=1}^R\hat{c}_j^{(m,r)}-c_j\right)^2
$$ was computed for $m=1,\dots,M$ and plotted in fig.2. We notice
that $e_R(m)\ll e_0(m)$ for almost all $m$ and it is much less
dispersed around its mean. Therefore the estimates of
$(\xi_j,c_j,j=1,\dots,p)$ obtained by averaging over the $R$
pseudosamples are better than those obtained by the original
samples. Finally we notice that in this simulation we used a
variance $\tilde{\sigma}^2$ much larger than the one used to produce
the results in Table 1. This large value gives the best mean squared
error over all the five parameters but not necessarily the best
reconstruction of each single parameter, as we looked for in the
previous simulation.

\begin{table}[tbh]
\begin{center}
\begin{tabular}{|c|c|c|c|c|}
\hline \hline
&$p$&$bias(\hat{p})$&$s.d.(\hat{p})$&$MSE(\hat{p})$\\
\hline
$ $&5& 0.0500&    1.0000 &   1.0025\\
\hline\hline
&$\xi_j$&$bias(\hat{\xi}_j)$&$s.d.\hat{\xi}_j$&$MSE(\hat{\xi}_j)$\\
\hline
$j=1$& -0.2796 - 0.8606i& -0.0006 + 0.0004i&   0.0230& 0.0005\\
\hline
$j=2$&-0.1782 - 0.9344i& -0.0005 - 0.0004i &  0.0125& 0.0002\\
\hline
$j=3$&   0.3090 + 0.9510i &  0.0057 - 0.0009i &  0.0171& 0.0003 \\
\hline
$j=4$&   0.2487 + 0.9685i & -0.0005 + 0.0024i &  0.0145 &0.0002 \\
\hline
$j=5$& -0.4354 + 0.5993i & -0.0054 + 0.0018i &  0.0290& 0.0009 \\
\hline \hline
&$c_j$&$bias(\hat{c}_j)$&$s.d.(\hat{c}_j)$&$MSE(\hat{c}_j)$\\
\hline
$j=1$&6.0000 &  0.1545  &  1.7154  &  2.9663\\
\hline
$j=2$&3.0000 &  -0.1617 &   1.2865 &   1.6812\\
\hline
$j=3$&1.0000 &  -0.1037  &  0.3295  &  0.1193\\
\hline
$j=4$&1.0000&  -0.0981 &   0.3193  &  0.1116\\
\hline
$j=5$& 20.0000& -0.1759 &   2.5101 &   6.3317\\
\hline
\end{tabular}
\end{center}
\caption{Statistics of the parameters $\hat{p}$,
$\hat{\xi}_j,j=1,\dots,p$ and $\hat{c}_j,j=1,\dots,p$ } \label{tb1}
\end{table}

\begin{figure}
\begin{center}
\fbox{\epsfig{file=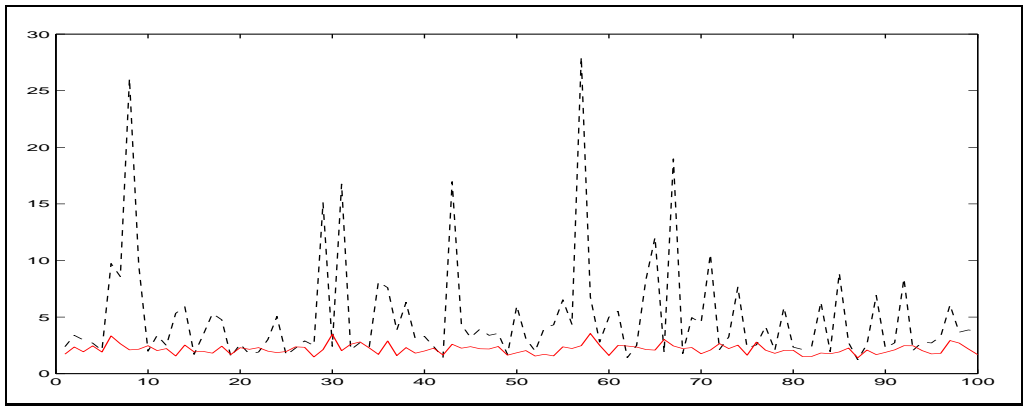,height=4.8cm,width=13cm}}
\end{center}
\caption{MSE of the standard estimator of the parameters
$(\xi_j,c_j),j=1,\dots,p$ (dashed); MSE of the averaged estimator
(solid)}\label{fig2}
\end{figure}

\end{document}